\documentclass[11pt]{article}
\usepackage{amsmath,amsthm,amsfonts,amssymb,url,graphicx,xcolor,authblk,longtable}
\usepackage{tikz}
\usepackage[title]{appendix}
\usetikzlibrary{calc}
\usetikzlibrary{arrows,decorations.markings}
\usepackage{hyperref}
\hypersetup{
    colorlinks=true,
    linkcolor=blue,
    filecolor=blue,
    urlcolor=blue,
    citecolor=blue,
}

\definecolor{red}{RGB}{220,20,60}% crimson

\newcommand{\MTS}{\ensuremath{\mathsf{MTS}}} 
\newcommand{\TTS}{\ensuremath{\mathsf{TTS}}}

\makeatletter
\g@addto@macro\bfseries{\boldmath}
\makeatother
\newcommand{\M}[3]{\ensuremath{\mathcal{M}_{#1}#2#3}}
\newcommand{\T}[2]{\ensuremath{\mathcal{T}_{#1}#2}}

\begin{document}
\title{\bf Good sequencings for\\small Mendelsohn  triple systems}
\author[1]{Donald L.\ Kreher}
\author[2]{Douglas R.\ Stinson%
\thanks{D.R.\ Stinson's research is supported by  NSERC discovery grant RGPIN-03882.}}
\author[3]{Shannon Veitch}
\affil[1]{Department of Mathematical Sciences, 
Michigan Technological University 
Houghton, MI 49931,  
U.S.A.}
\affil[2]{David R.\ Cheriton School of Computer Science, University of Waterloo,
Waterloo, Ontario, N2L 3G1, Canada}
\affil[3]{Department of Combinatorics and Optimization, University of Waterloo,
Waterloo, Ontario, N2L 3G1, Canada}

\maketitle

\begin{abstract}
A \emph{Mendelsohn triple system} of order $v$ (or $\MTS(v)$) is a 
decomposition of the complete graph into directed 3-cyles. We denote 
the directed 3-cycle with edges $(x,y)$, $(y,z)$ and $(z,x)$ by 
$(x,y,z)$, $(y,z,x)$ or $(z,x,y)$. 

An \emph{$\ell$-good sequencing}  of a \MTS$(v)$  
is a permutation of the points of the design, 
say $[x_1 \; \cdots \; x_v]$, such that, for every 
triple $(x,y,z)$ in the design,
it is \emph{not} the case that 
$x = x_i$, $y = x_j$ and $z = x_k$ 
   with $i < j < k$  and $k-i+1 \leq \ell$;
or with $j < k < i$  and $i-j+1 \leq \ell$;
or with $k < i < j$  and $j-k+1 \leq \ell$.

In this report we provide a maximum $\ell$-good sequencing for each 
\MTS$(v)$, $v \leq 10$.
\end{abstract}

\section{Introduction}
It is in~\cite{KSV} that  
\emph{$\ell$-good sequencings} of Mendelsohn triple systems are first 
introduced and this paper is where we direct the interested reader to 
find definitions, fundamental results and constructions. In particular,
it is well known that an \MTS$(v)$ exists if and only if 
$v \equiv 0 \text{ or } 1 \pmod{3}$, $v \neq 6$. 
In~\cite{KSV}, it is shown that, if an $\MTS(v)$ has
an $\ell$-good sequencing, then $\ell \leq \lfloor \frac{v-1}{2} \rfloor$.

In \hyperref[summary]{Table~\ref{summary}}, we report, for $v\leq 10$, the number 
of nonisomorphic $\MTS(v)$ that have an $\ell$-good sequencing. 
If each triple $(x,y,z)$ in a Mendelsohn triple system is replaced with the 3-element subset $\{x,y,z\}$, a twofold triple 
system of order $v$ (i.e., a $\TTS(v)$) is obtained. 
This is a decomposition of $2K_v$,
the double of the complete graph, into triangles.
Alternatively, starting with a $\TTS(v)$, it may be possible to orient each of the triangles to obtain a $\MTS(v)$. Such $\TTS(v)$
are said to be orientable.
The numbers of nonisomorphic $\TTS(v)$, orientable $\TTS(v)$ 
and nonisomorphic $\MTS(v)$ were previously reported 
in~\cite{DG}. We confirm their results.

\begin{table}[h]
\caption{Summary of results}\label{summary}
\begin{center}
\begin{tabular}{@{}*{7}{r}@{}}
\hline
  & 
  Nonisomorphic& 
  Orientable&
  Nonisomorphic& 
  \multicolumn{3}{@{}c@{}}{$\ell$-good sequencing}\\
  $v$&
  $\TTS(v)$&
  $\TTS(v)$&
  $\MTS(v)$&
  $\ell=3$ & $\ell=4$ & $\ell=5$\\
\hline
  3& 1&   1&  1&  0&  0&0 \\
  4& 1&   1&  1&  0&  0&0 \\
  6& 1&   0&  0&  0&  0&0 \\
  7& 4&   3&  3&  3&  0&0 \\
  9&36&  16& 18& 18&  3&0 \\
 10&960&134&143&143&138&0 \\
 \hline
\end{tabular}
\end{center}
\end{table}

In the following sections, we denote the $i$-th orientable 
$\TTS(v)$ by $\T{v}{i}$.  The $j$-th $\MTS(v)$ 
whose underlying $\TTS(v)$ is $\T{v}{i}$,
is denoted by $\M{v}{i}{.j}$.

\newpage
\section{Orientable \texorpdfstring{$\TTS(7)$s}{TTS(7)s}}\label{v7}
There are 3 $\TTS(7)$ that are orientable. They are listed in \hyperref[TabI7]{Section~\ref{TabI7}}.
The $\MTS(7)$ that can be obtained by orienting the $\TTS(7)$ are listed in 
\hyperref[TabII7]{Section~\ref{TabII7}}. 
\subsection{Orientable \texorpdfstring{$\TTS(7)$}{TTS(7)}s.}\label{TabI7}

\vbox{\begin{description}
\item[\T{7}{1}:]$\begin{array}{@{}*{6}{l@{\hspace{4pt}}}l@{}}
 \{3,4,5\}& \{2,5,6\}& \{2,3,6\}& \{2,3,4\}& \{1,4,6\}& \{1,4,5\}& \{1,3,6\}\\
 \{1,3,6\}& \{1,2,5\}& \{0,5,6\}& \{0,4,6\}& \{0,3,5\}& \{0,2,4\}& \{0,1,3\}
\end{array}$
\end{description}}
\vbox{\begin{description}
\item[\T{7}{2}:]$\begin{array}{@{}*{6}{l@{\hspace{4pt}}}l@{}}
 \{2,4,5\}& \{2,4,5\}& \{2,3,6\}& \{2,3,6\}& \{1,4,6\}& \{1,4,6\}& \{1,3,5\}\\
 \{1,3,5\}& \{1,3,5\}& \{0,5,6\}& \{0,5,6\}& \{0,3,4\}& \{0,3,4\}& \{0,1,2\}
\end{array}$
\end{description}}
\vbox{\begin{description}
\item[\T{7}{3}:]$\begin{array}{@{}*{6}{l@{\hspace{4pt}}}l@{}}
 \{2,4,5\}& \{2,4,5\}& \{2,3,6\}& \{2,3,6\}& \{1,5,6\}& \{1,4,6\}& \{1,3,5\}\\
 \{1,3,5\}& \{1,3,4\}& \{0,5,6\}& \{0,4,6\}& \{0,3,5\}& \{0,3,4\}& \{0,1,2\}
\end{array}$
\end{description}}
\subsection{The \texorpdfstring{$\MTS(7)$}{MTS(7)} and the good sequencings found.}\label{TabII7}

\vbox{\begin{description}
\item[\M{7}{1}{.1}:]$\begin{array}{@{}*{6}{l@{\hspace{4pt}}}l@{}}
 (0,1,2)& (0,3,1)& (0,2,4)& (0,5,3)& (0,4,6)& (0,6,5)& (1,5,2)\\
 (1,3,6)& (1,4,5)& (1,6,4)& (2,3,4)& (2,6,3)& (2,5,6)& (3,5,4)
\end{array}$

	Lexicographic least 3-good sequencing : 0132456

	Number of 3-good sequencing found: 336
\end{description}}
\vbox{\begin{description}
\item[\M{7}{2}{.1}:]$\begin{array}{@{}*{6}{l@{\hspace{4pt}}}l@{}}
 (0,1,2)& (0,2,1)& (0,3,4)& (0,4,3)& (0,5,6)& (0,6,5)& (1,3,5)\\
 (1,5,3)& (1,4,6)& (1,6,4)& (2,3,6)& (2,6,3)& (2,4,5)& (2,5,4)
\end{array}$

	Lexicographic least 3-good sequencing : 0132564

	Number of 3-good sequencing found: 336
\end{description}}
\vbox{\begin{description}
\item[\M{7}{3}{.1}:]$\begin{array}{@{}*{6}{l@{\hspace{4pt}}}l@{}}
 (0,1,2)& (0,2,1)& (0,4,3)& (0,3,5)& (0,6,4)& (0,5,6)& (1,3,4)\\
 (1,5,3)& (1,4,6)& (1,6,5)& (2,3,6)& (2,6,3)& (2,4,5)& (2,5,4)
\end{array}$

	Lexicographic least 3-good sequencing : 0132465

	Number of 3-good sequencing found: 336
\end{description}}

\newpage
\section{Orientable \texorpdfstring{$\TTS(9)$}{TTS(9)}s}\label{v9}
There are 16 $\TTS(9)$ that are orientable. They are listed in \hyperref[TabI9]{Section~\ref{TabI9}}.
The $\MTS(9)$ that can be obtained by orienting the $\TTS(9)$ are listed in 
\hyperref[TabII9]{Section~\ref{TabII9}}.

\subsection{Orientable \texorpdfstring{$\TTS(9)$}{TTS(9)}.}\label{TabI9}

\vbox{\begin{description}
\item[\T{9}{1}:]	
$\begin{array}{@{}*{7}{l@{\hspace{4pt}}}l@{}}
 \{6,7,8\}& \{6,7,8\}& \{3,5,6\}& \{3,4,7\}& \{3,4,5\}& \{2,5,8\}& \{2,5,7\}& \{2,4,6\}\\
 \{2,4,6\}& \{2,4,6\}& \{2,3,8\}& \{1,5,7\}& \{1,4,8\}& \{1,4,5\}& \{1,3,8\}& \{1,3,6\}\\
 \{1,3,6\}& \{1,2,7\}& \{0,5,8\}& \{0,5,6\}& \{0,4,8\}& \{0,4,7\}& \{0,3,7\}& \{0,2,3\}
\end{array}$
\end{description}}
\vbox{\begin{description}
\item[\T{9}{2}:]	
$\begin{array}{@{}*{7}{l@{\hspace{4pt}}}l@{}}
 \{6,7,8\}& \{6,7,8\}& \{4,5,6\}& \{3,5,8\}& \{3,4,7\}& \{3,4,7\}& \{2,5,8\}& \{2,4,8\}\\
 \{2,4,8\}& \{2,3,6\}& \{2,3,5\}& \{1,5,7\}& \{1,4,8\}& \{1,4,5\}& \{1,3,6\}& \{1,2,7\}\\
 \{1,2,7\}& \{1,2,6\}& \{0,5,7\}& \{0,5,6\}& \{0,4,6\}& \{0,3,8\}& \{0,2,7\}& \{0,2,4\}
\end{array}$
\end{description}}
\vbox{\begin{description}
\item[\T{9}{3}:]	
$\begin{array}{@{}*{7}{l@{\hspace{4pt}}}l@{}}
 \{6,7,8\}& \{6,7,8\}& \{3,5,6\}& \{3,4,7\}& \{3,4,5\}& \{2,5,8\}& \{2,4,8\}& \{2,4,5\}\\
 \{2,4,5\}& \{2,3,7\}& \{2,3,6\}& \{1,5,7\}& \{1,5,6\}& \{1,4,8\}& \{1,4,6\}& \{1,3,8\}\\
 \{1,3,8\}& \{1,2,7\}& \{0,5,8\}& \{0,5,7\}& \{0,4,7\}& \{0,4,6\}& \{0,3,8\}& \{0,2,6\}
\end{array}$
\end{description}}
\vbox{\begin{description}
\item[\T{9}{4}:]	
$\begin{array}{@{}*{7}{l@{\hspace{4pt}}}l@{}}
 \{6,7,8\}& \{6,7,8\}& \{4,5,6\}& \{3,5,7\}& \{3,4,8\}& \{3,4,8\}& \{2,5,8\}& \{2,5,8\}\\
 \{2,5,8\}& \{2,4,7\}& \{2,3,6\}& \{1,5,7\}& \{1,4,7\}& \{1,3,6\}& \{1,3,5\}& \{1,2,6\}\\
 \{1,2,6\}& \{1,2,4\}& \{0,5,6\}& \{0,4,6\}& \{0,4,5\}& \{0,3,7\}& \{0,2,7\}& \{0,2,3\}
\end{array}$
\end{description}}
\vbox{\begin{description}
\item[\T{9}{5}:]	
$\begin{array}{@{}*{7}{l@{\hspace{4pt}}}l@{}}
 \{6,7,8\}& \{6,7,8\}& \{4,5,6\}& \{3,5,7\}& \{3,4,8\}& \{2,5,8\}& \{2,4,8\}& \{2,4,7\}\\
 \{2,4,7\}& \{2,3,6\}& \{2,3,5\}& \{1,5,8\}& \{1,5,7\}& \{1,4,7\}& \{1,3,6\}& \{1,3,4\}\\
 \{1,3,4\}& \{1,2,6\}& \{0,5,6\}& \{0,4,6\}& \{0,4,5\}& \{0,3,8\}& \{0,3,7\}& \{0,2,7\}
\end{array}$
\end{description}}
\vbox{\begin{description}
\item[\T{9}{6}:]	
$\begin{array}{@{}*{7}{l@{\hspace{4pt}}}l@{}}
 \{6,7,8\}& \{6,7,8\}& \{3,5,6\}& \{3,4,7\}& \{3,4,5\}& \{2,5,8\}& \{2,5,7\}& \{2,4,8\}\\
 \{2,4,8\}& \{2,4,6\}& \{2,3,6\}& \{1,5,7\}& \{1,4,6\}& \{1,4,5\}& \{1,3,8\}& \{1,3,8\}\\
 \{1,3,8\}& \{1,2,7\}& \{0,5,8\}& \{0,5,6\}& \{0,4,8\}& \{0,4,7\}& \{0,3,7\}& \{0,2,3\}
\end{array}$
\end{description}}
\vbox{\begin{description}
\item[\T{9}{7}:]	
$\begin{array}{@{}*{7}{l@{\hspace{4pt}}}l@{}}
 \{6,7,8\}& \{6,7,8\}& \{3,4,5\}& \{3,4,5\}& \{2,5,7\}& \{2,5,6\}& \{2,4,8\}& \{2,4,6\}\\
 \{2,4,6\}& \{2,3,8\}& \{2,3,7\}& \{1,5,8\}& \{1,5,6\}& \{1,4,8\}& \{1,4,7\}& \{1,3,7\}\\
 \{1,3,7\}& \{1,3,6\}& \{0,5,8\}& \{0,5,7\}& \{0,4,7\}& \{0,4,6\}& \{0,3,8\}& \{0,3,6\}
\end{array}$
\end{description}}
\vbox{\begin{description}
\item[\T{9}{8}:]	
$\begin{array}{@{}*{7}{l@{\hspace{4pt}}}l@{}}
 \{6,7,8\}& \{6,7,8\}& \{3,5,6\}& \{3,4,7\}& \{3,4,5\}& \{2,5,7\}& \{2,4,6\}& \{2,4,5\}\\
 \{2,4,5\}& \{2,3,8\}& \{2,3,8\}& \{1,5,8\}& \{1,5,8\}& \{1,4,7\}& \{1,4,6\}& \{1,3,6\}\\
 \{1,3,6\}& \{1,2,7\}& \{0,5,7\}& \{0,5,6\}& \{0,4,8\}& \{0,4,8\}& \{0,3,7\}& \{0,2,6\}
\end{array}$
\end{description}}
\vbox{\begin{description}
\item[\T{9}{9}:]	
$\begin{array}{@{}*{7}{l@{\hspace{4pt}}}l@{}}
 \{6,7,8\}& \{6,7,8\}& \{3,5,6\}& \{3,4,7\}& \{3,4,5\}& \{2,5,7\}& \{2,4,6\}& \{2,4,5\}\\
 \{2,4,5\}& \{2,3,8\}& \{2,3,8\}& \{1,5,8\}& \{1,5,8\}& \{1,4,7\}& \{1,4,6\}& \{1,3,7\}\\
 \{1,3,7\}& \{1,2,6\}& \{0,5,7\}& \{0,5,6\}& \{0,4,8\}& \{0,4,8\}& \{0,3,6\}& \{0,2,7\}
\end{array}$
\end{description}}
\vbox{\begin{description}
\item[\T{9}{10}:]	
$\begin{array}{@{}*{7}{l@{\hspace{4pt}}}l@{}}
 \{6,7,8\}& \{6,7,8\}& \{3,4,5\}& \{3,4,5\}& \{2,5,6\}& \{2,5,6\}& \{2,4,7\}& \{2,4,7\}\\
 \{2,4,7\}& \{2,3,8\}& \{2,3,8\}& \{1,5,8\}& \{1,5,7\}& \{1,4,8\}& \{1,4,6\}& \{1,3,7\}\\
 \{1,3,7\}& \{1,3,6\}& \{0,5,8\}& \{0,5,7\}& \{0,4,8\}& \{0,4,6\}& \{0,3,7\}& \{0,3,6\}
\end{array}$
\end{description}}
\vbox{\begin{description}
\item[\T{9}{11}:]	
$\begin{array}{@{}*{7}{l@{\hspace{4pt}}}l@{}}
 \{6,7,8\}& \{6,7,8\}& \{3,4,5\}& \{3,4,5\}& \{2,5,8\}& \{2,5,8\}& \{2,4,6\}& \{2,4,6\}\\
 \{2,4,6\}& \{2,3,7\}& \{2,3,7\}& \{1,5,6\}& \{1,5,6\}& \{1,4,7\}& \{1,4,7\}& \{1,3,8\}\\
 \{1,3,8\}& \{1,3,8\}& \{0,5,7\}& \{0,5,7\}& \{0,4,8\}& \{0,4,8\}& \{0,3,6\}& \{0,3,6\}
\end{array}$
\end{description}}
\vbox{\begin{description}
\item[\T{9}{12}:]	
$\begin{array}{@{}*{7}{l@{\hspace{4pt}}}l@{}}
 \{4,5,7\}& \{4,5,6\}& \{3,6,7\}& \{3,5,8\}& \{3,4,8\}& \{3,4,6\}& \{2,6,8\}& \{2,5,8\}\\
 \{2,5,8\}& \{2,4,7\}& \{2,3,7\}& \{1,7,8\}& \{1,6,7\}& \{1,4,8\}& \{1,3,5\}& \{1,2,6\}\\
 \{1,2,6\}& \{1,2,5\}& \{0,7,8\}& \{0,6,8\}& \{0,5,7\}& \{0,5,6\}& \{0,2,4\}& \{0,2,3\}
\end{array}$
\end{description}}
\vbox{\begin{description}
\item[\T{9}{13}:]	
$\begin{array}{@{}*{7}{l@{\hspace{4pt}}}l@{}}
 \{3,6,7\}& \{3,5,8\}& \{3,5,6\}& \{3,4,8\}& \{3,4,7\}& \{2,6,7\}& \{2,5,8\}& \{2,5,6\}\\
 \{2,5,6\}& \{2,4,8\}& \{2,4,7\}& \{1,7,8\}& \{1,6,8\}& \{1,5,7\}& \{1,4,6\}& \{1,4,5\}\\
 \{1,4,5\}& \{1,2,3\}& \{0,7,8\}& \{0,6,8\}& \{0,5,7\}& \{0,4,6\}& \{0,4,5\}& \{0,2,3\}
\end{array}$
\end{description}}
\vbox{\begin{description}
\item[\T{9}{14}:]	
$\begin{array}{@{}*{7}{l@{\hspace{4pt}}}l@{}}
 \{4,5,6\}& \{3,6,7\}& \{3,5,8\}& \{3,4,8\}& \{3,4,7\}& \{2,6,8\}& \{2,5,8\}& \{2,5,7\}\\
 \{2,5,7\}& \{2,4,7\}& \{2,3,6\}& \{1,7,8\}& \{1,5,7\}& \{1,5,6\}& \{1,4,8\}& \{1,4,6\}\\
 \{1,4,6\}& \{1,2,3\}& \{0,7,8\}& \{0,6,8\}& \{0,6,7\}& \{0,4,5\}& \{0,3,5\}& \{0,2,4\}
\end{array}$
\end{description}}
\vbox{\begin{description}
\item[\T{9}{15}:]	
$\begin{array}{@{}*{7}{l@{\hspace{4pt}}}l@{}}
 \{4,5,7\}& \{4,5,6\}& \{3,6,7\}& \{3,5,8\}& \{3,4,8\}& \{2,6,8\}& \{2,5,7\}& \{2,5,6\}\\
 \{2,5,6\}& \{2,4,8\}& \{2,3,7\}& \{1,7,8\}& \{1,6,7\}& \{1,5,8\}& \{1,4,6\}& \{1,3,4\}\\
 \{1,3,4\}& \{1,2,3\}& \{0,7,8\}& \{0,6,8\}& \{0,4,7\}& \{0,3,6\}& \{0,3,5\}& \{0,2,4\}
\end{array}$
\end{description}}
\vbox{\begin{description}
\item[\T{9}{16}:]	
$\begin{array}{@{}*{7}{l@{\hspace{4pt}}}l@{}}
 \{4,5,6\}& \{3,6,7\}& \{3,5,8\}& \{3,4,8\}& \{3,4,7\}& \{2,6,8\}& \{2,5,7\}& \{2,5,6\}\\
 \{2,5,6\}& \{2,4,8\}& \{2,4,7\}& \{1,7,8\}& \{1,6,7\}& \{1,5,8\}& \{1,4,5\}& \{1,3,6\}\\
 \{1,3,6\}& \{1,2,3\}& \{0,7,8\}& \{0,6,8\}& \{0,5,7\}& \{0,4,6\}& \{0,3,5\}& \{0,2,3\}
\end{array}$
\end{description}}
\subsection{The \texorpdfstring{$\TTS(9)$}{TTS(9)} and the good sequencings found.}\label{TabII9}

\vbox{\begin{description}
\item[\M{9}{1}{.1}:]	
$\begin{array}{@{}*{7}{l@{\hspace{4pt}}}l@{}}
 (0,2,1)& (0,1,6)& (0,3,2)& (0,7,3)& (0,4,7)& (0,8,4)& (0,6,5)& (0,5,8)\\
 (1,2,7)& (1,3,6)& (1,8,3)& (1,5,4)& (1,4,8)& (1,7,5)& (2,3,8)& (2,4,6)\\
 (2,6,4)& (2,5,7)& (2,8,5)& (3,4,5)& (3,7,4)& (3,5,6)& (6,7,8)& (6,8,7)
\end{array}$

	Lexicographic least 3-good sequencing : 012346578

	Number of 3-good sequencing found: 60156

	Lexicographic least 4-good sequencing : 023471856

	Number of 4-good sequencing found: 18
\end{description}}
\vbox{\begin{description}
\item[\M{9}{2}{.1}:]	
$\begin{array}{@{}*{7}{l@{\hspace{4pt}}}l@{}}
 (0,3,1)& (0,1,8)& (0,2,4)& (0,7,2)& (0,8,3)& (0,4,6)& (0,6,5)& (0,5,7)\\
 (1,6,2)& (1,2,7)& (1,3,6)& (1,5,4)& (1,4,8)& (1,7,5)& (2,3,5)& (2,6,3)\\
 (2,8,4)& (2,5,8)& (3,4,7)& (3,7,4)& (3,8,5)& (4,5,6)& (6,7,8)& (6,8,7)
\end{array}$

	Lexicographic least 3-good sequencing : 012346587

	Number of 3-good sequencing found: 60840

	Number of 4-good sequencing found: 0
\end{description}}
\vbox{\begin{description}
\item[\M{9}{3}{.1}:]	
$\begin{array}{@{}*{7}{l@{\hspace{4pt}}}l@{}}
 (0,2,1)& (0,1,3)& (0,6,2)& (0,3,8)& (0,4,6)& (0,7,4)& (0,5,7)& (0,8,5)\\
 (1,2,7)& (1,8,3)& (1,6,4)& (1,4,8)& (1,5,6)& (1,7,5)& (2,6,3)& (2,3,7)\\
 (2,4,5)& (2,8,4)& (2,5,8)& (3,5,4)& (3,4,7)& (3,6,5)& (6,7,8)& (6,8,7)
\end{array}$

	Lexicographic least 3-good sequencing : 012346578

	Number of 3-good sequencing found: 60696

	Lexicographic least 4-good sequencing : 047563812

	Number of 4-good sequencing found: 36
\end{description}}

\vbox{\begin{description}
\item[\M{9}{4}{.1}:]	
$\begin{array}{@{}*{7}{l@{\hspace{4pt}}}l@{}}
 (0,1,8)& (0,8,1)& (0,3,2)& (0,2,7)& (0,7,3)& (0,5,4)& (0,4,6)& (0,6,5)\\
 (1,4,2)& (1,2,6)& (1,3,5)& (1,6,3)& (1,7,4)& (1,5,7)& (2,3,6)& (2,4,7)\\
 (2,5,8)& (2,8,5)& (3,4,8)& (3,8,4)& (3,7,5)& (4,5,6)& (6,7,8)& (6,8,7)
\end{array}$

	Lexicographic least 3-good sequencing : 012346587

	Number of 3-good sequencing found: 60480

	Number of 4-good sequencing found: 0
\end{description}}
\vbox{\begin{description}
\item[\M{9}{4}{.2}:]	
$\begin{array}{@{}*{7}{l@{\hspace{4pt}}}l@{}}
 (0,1,8)& (0,8,1)& (0,3,2)& (0,2,7)& (0,7,3)& (0,4,5)& (0,6,4)& (0,5,6)\\
 (1,4,2)& (1,2,6)& (1,3,5)& (1,6,3)& (1,7,4)& (1,5,7)& (2,3,6)& (2,4,7)\\
 (2,5,8)& (2,8,5)& (3,4,8)& (3,8,4)& (3,7,5)& (4,6,5)& (6,7,8)& (6,8,7)
\end{array}$

	Lexicographic least 3-good sequencing : 012346857

	Number of 3-good sequencing found: 59616

	Number of 4-good sequencing found: 0
\end{description}}
\vbox{\begin{description}
\item[\M{9}{5}{.1}:]	
$\begin{array}{@{}*{7}{l@{\hspace{4pt}}}l@{}}
 (0,2,1)& (0,1,8)& (0,7,2)& (0,3,7)& (0,8,3)& (0,5,4)& (0,4,6)& (0,6,5)\\
 (1,2,6)& (1,3,4)& (1,6,3)& (1,4,7)& (1,7,5)& (1,5,8)& (2,5,3)& (2,3,6)\\
 (2,7,4)& (2,4,8)& (2,8,5)& (3,8,4)& (3,5,7)& (4,5,6)& (6,7,8)& (6,8,7)
\end{array}$

	Lexicographic least 3-good sequencing : 012346587

	Number of 3-good sequencing found: 61506

	Number of 4-good sequencing found: 0
\end{description}}
\vbox{\begin{description}
\item[\M{9}{5}{.2}:]	
$\begin{array}{@{}*{7}{l@{\hspace{4pt}}}l@{}}
 (0,2,1)& (0,1,8)& (0,7,2)& (0,3,7)& (0,8,3)& (0,4,5)& (0,6,4)& (0,5,6)\\
 (1,2,6)& (1,3,4)& (1,6,3)& (1,4,7)& (1,7,5)& (1,5,8)& (2,5,3)& (2,3,6)\\
 (2,7,4)& (2,4,8)& (2,8,5)& (3,8,4)& (3,5,7)& (4,6,5)& (6,7,8)& (6,8,7)
\end{array}$

	Lexicographic least 3-good sequencing : 012346857

	Number of 3-good sequencing found: 61290

	Number of 4-good sequencing found: 0
\end{description}}
\vbox{\begin{description}
\item[\M{9}{6}{.1}:]	
$\begin{array}{@{}*{7}{l@{\hspace{4pt}}}l@{}}
 (0,2,1)& (0,1,6)& (0,3,2)& (0,7,3)& (0,4,7)& (0,8,4)& (0,6,5)& (0,5,8)\\
 (1,2,7)& (1,3,8)& (1,8,3)& (1,5,4)& (1,4,6)& (1,7,5)& (2,3,6)& (2,6,4)\\
 (2,4,8)& (2,5,7)& (2,8,5)& (3,4,5)& (3,7,4)& (3,5,6)& (6,7,8)& (6,8,7)
\end{array}$

	Lexicographic least 3-good sequencing : 012346578

	Number of 3-good sequencing found: 60300

	Number of 4-good sequencing found: 0
\end{description}}
\vbox{\begin{description}
\item[\M{9}{7}{.1}:]	
$\begin{array}{@{}*{7}{l@{\hspace{4pt}}}l@{}}
 (0,1,2)& (0,2,1)& (0,6,3)& (0,3,8)& (0,4,6)& (0,7,4)& (0,5,7)& (0,8,5)\\
 (1,3,6)& (1,7,3)& (1,4,7)& (1,8,4)& (1,6,5)& (1,5,8)& (2,3,7)& (2,8,3)\\
 (2,6,4)& (2,4,8)& (2,5,6)& (2,7,5)& (3,4,5)& (3,5,4)& (6,7,8)& (6,8,7)
\end{array}$

	Lexicographic least 3-good sequencing : 013246578

	Number of 3-good sequencing found: 58644

	Lexicographic least 4-good sequencing : 031485726

	Number of 4-good sequencing found: 324
\end{description}}
\vbox{\begin{description}
\item[\M{9}{8}{.1}:]	
$\begin{array}{@{}*{7}{l@{\hspace{4pt}}}l@{}}
 (0,2,1)& (0,1,3)& (0,6,2)& (0,3,7)& (0,4,8)& (0,8,4)& (0,5,6)& (0,7,5)\\
 (1,2,7)& (1,6,3)& (1,4,6)& (1,7,4)& (1,5,8)& (1,8,5)& (2,3,8)& (2,8,3)\\
 (2,4,5)& (2,6,4)& (2,5,7)& (3,5,4)& (3,4,7)& (3,6,5)& (6,7,8)& (6,8,7)
\end{array}$

	Lexicographic least 3-good sequencing : 012346578

	Number of 3-good sequencing found: 58176

	Number of 4-good sequencing found: 0
\end{description}}
\vbox{\begin{description}
\item[\M{9}{9}{.1}:]	
$\begin{array}{@{}*{7}{l@{\hspace{4pt}}}l@{}}
 (0,2,1)& (0,1,3)& (0,7,2)& (0,3,6)& (0,4,8)& (0,8,4)& (0,6,5)& (0,5,7)\\
 (1,2,6)& (1,7,3)& (1,6,4)& (1,4,7)& (1,5,8)& (1,8,5)& (2,3,8)& (2,8,3)\\
 (2,5,4)& (2,4,6)& (2,7,5)& (3,4,5)& (3,7,4)& (3,5,6)& (6,7,8)& (6,8,7)
\end{array}$

	Lexicographic least 3-good sequencing : 012346578

	Number of 3-good sequencing found: 59472

	Number of 4-good sequencing found: 0
\end{description}}
\vbox{\begin{description}
\item[\M{9}{10}{.1}:]	
$\begin{array}{@{}*{7}{l@{\hspace{4pt}}}l@{}}
 (0,1,2)& (0,2,1)& (0,6,3)& (0,3,7)& (0,4,6)& (0,8,4)& (0,7,5)& (0,5,8)\\
 (1,3,6)& (1,7,3)& (1,6,4)& (1,4,8)& (1,5,7)& (1,8,5)& (2,3,8)& (2,8,3)\\
 (2,4,7)& (2,7,4)& (2,5,6)& (2,6,5)& (3,4,5)& (3,5,4)& (6,7,8)& (6,8,7)
\end{array}$

	Lexicographic least 3-good sequencing : 013246578

	Number of 3-good sequencing found: 56592

	Number of 4-good sequencing found: 0
\end{description}}
\vbox{\begin{description}
\item[\M{9}{11}{.1}:]	
$\begin{array}{@{}*{7}{l@{\hspace{4pt}}}l@{}}
 (0,1,2)& (0,2,1)& (0,3,6)& (0,6,3)& (0,4,8)& (0,8,4)& (0,5,7)& (0,7,5)\\
 (1,3,8)& (1,8,3)& (1,4,7)& (1,7,4)& (1,5,6)& (1,6,5)& (2,3,7)& (2,7,3)\\
 (2,4,6)& (2,6,4)& (2,5,8)& (2,8,5)& (3,4,5)& (3,5,4)& (6,7,8)& (6,8,7)
\end{array}$

	Lexicographic least 3-good sequencing : 013247568

	Number of 3-good sequencing found: 51840

	Number of 4-good sequencing found: 0
\end{description}}
\vbox{\begin{description}
\item[\M{9}{12}{.1}:]	
$\begin{array}{@{}*{7}{l@{\hspace{4pt}}}l@{}}
 (0,3,1)& (0,1,4)& (0,2,3)& (0,4,2)& (0,5,6)& (0,7,5)& (0,6,8)& (0,8,7)\\
 (1,5,2)& (1,2,6)& (1,3,5)& (1,8,4)& (1,6,7)& (1,7,8)& (2,7,3)& (2,4,7)\\
 (2,5,8)& (2,8,6)& (3,6,4)& (3,4,8)& (3,8,5)& (3,7,6)& (4,6,5)& (4,5,7)
\end{array}$

	Lexicographic least 3-good sequencing : 012345678

	Number of 3-good sequencing found: 61128

	Number of 4-good sequencing found: 0
\end{description}}
\vbox{\begin{description}
\item[\M{9}{13}{.1}:]	
$\begin{array}{@{}*{7}{l@{\hspace{4pt}}}l@{}}
 (0,2,1)& (0,1,3)& (0,3,2)& (0,5,4)& (0,4,6)& (0,7,5)& (0,6,8)& (0,8,7)\\
 (1,2,3)& (1,4,5)& (1,6,4)& (1,5,7)& (1,8,6)& (1,7,8)& (2,7,4)& (2,4,8)\\
 (2,5,6)& (2,8,5)& (2,6,7)& (3,4,7)& (3,8,4)& (3,6,5)& (3,5,8)& (3,7,6)
\end{array}$

	Lexicographic least 3-good sequencing : 012435678

	Number of 3-good sequencing found: 61920

	Number of 4-good sequencing found: 0
\end{description}}
\vbox{\begin{description}
\item[\M{9}{14}{.1}:]	
$\begin{array}{@{}*{7}{l@{\hspace{4pt}}}l@{}}
 (0,2,1)& (0,1,3)& (0,4,2)& (0,3,5)& (0,5,4)& (0,7,6)& (0,6,8)& (0,8,7)\\
 (1,2,3)& (1,4,6)& (1,8,4)& (1,6,5)& (1,5,7)& (1,7,8)& (2,6,3)& (2,4,7)\\
 (2,7,5)& (2,5,8)& (2,8,6)& (3,7,4)& (3,4,8)& (3,8,5)& (3,6,7)& (4,5,6)
\end{array}$

	Lexicographic least 3-good sequencing : 012435678

	Number of 3-good sequencing found: 61614

	Number of 4-good sequencing found: 0
\end{description}}
\vbox{\begin{description}
\item[\M{9}{15}{.1}:]	
$\begin{array}{@{}*{7}{l@{\hspace{4pt}}}l@{}}
 (0,2,1)& (0,1,5)& (0,4,2)& (0,5,3)& (0,3,6)& (0,7,4)& (0,6,8)& (0,8,7)\\
 (1,2,3)& (1,3,4)& (1,4,6)& (1,8,5)& (1,6,7)& (1,7,8)& (2,7,3)& (2,4,8)\\
 (2,6,5)& (2,5,7)& (2,8,6)& (3,8,4)& (3,5,8)& (3,7,6)& (4,5,6)& (4,7,5)
\end{array}$

	Lexicographic least 3-good sequencing : 012435678

	Number of 3-good sequencing found: 60912

	Number of 4-good sequencing found: 0
\end{description}}
\vbox{\begin{description}
\item[\M{9}{16}{.1}:]	
$\begin{array}{@{}*{7}{l@{\hspace{4pt}}}l@{}}
 (0,2,1)& (0,1,4)& (0,3,2)& (0,5,3)& (0,4,6)& (0,7,5)& (0,6,8)& (0,8,7)\\
 (1,2,3)& (1,3,6)& (1,5,4)& (1,8,5)& (1,6,7)& (1,7,8)& (2,7,4)& (2,4,8)\\
 (2,6,5)& (2,5,7)& (2,8,6)& (3,4,7)& (3,8,4)& (3,5,8)& (3,7,6)& (4,5,6)
\end{array}$

	Lexicographic least 3-good sequencing : 012435678

	Number of 3-good sequencing found: 61632

	Number of 4-good sequencing found: 0
\end{description}}

\newpage
\section{Orientable \texorpdfstring{$\TTS(10)$}{TTS(10)}s}\label{v10}
There are 134 $\TTS(10)$, that are orientable. They are listed in \hyperref[TabI10]{Section~\ref{TabI10}}
and labeled $\T{10}{i}$, $i=1,2,\dots,134$. 
Following label $\T{10}{i}$ in brackets
is  the label of the same twofold triple system
 found in~\cite{GMRI} and \cite{GMRII}.
The $\MTS(10)$ that can be obtained by orienting the $\TTS(10)$ are listed in \hyperref[TabII10]{Section~\ref{TabII10}}.  The five $\MTS(10)$ that
do not have a 4-good sequencing are
%M116-1 M116-2 M118-1 M134-1 M134-2
$\M{10}{116}{.1}$, $\M{10}{116}{.2}$, $\M{10}{118}{.1}$, $\M{10}{134}{.1}$ and $\M{10}{134}{.2}$. If any triple is omitted from one of these five $\MTS(10)$, then the resulting ``partial'' 
$\MTS(10)$ having 29 triples turns out to have 
a 4-good sequencing. These 4-good sequencings are reported in
\hyperref[Tab29]{Section~\ref{Tab29}}.

\subsection{Orientable \texorpdfstring{$\TTS(10)$}{TTS(10)}.}\label{TabI10}

\vbox{\begin{description}
\item[\T{10}{1}:] [(1)]\\

$
\begin{array}{@{}*{5}{l@{\hspace{4pt}}}l@{}}
 \{1,2,3\}& \{2,5,7\}& \{3,7,9\}& \{1,3,4\}& \{2,5,6\}& \{0,3,6\}\\
 \{0,3,6\}& \{1,2,4\}& \{2,7,9\}& \{4,5,9\}& \{1,5,6\}& \{2,6,8\}\\
 \{2,6,8\}& \{4,6,8\}& \{1,6,7\}& \{0,2,9\}& \{0,4,7\}& \{1,5,7\}\\
 \{1,5,7\}& \{0,2,8\}& \{4,7,8\}& \{0,1,8\}& \{0,3,5\}& \{4,6,9\}\\
 \{4,6,9\}& \{1,8,9\}& \{3,6,9\}& \{0,4,5\}& \{0,1,9\}& \{3,7,8\}
\end{array}$
\end{description}}
\vbox{\begin{description}
\item[\T{10}{2}:] [(26)]\\

$\begin{array}{@{}*{5}{l@{\hspace{4pt}}}l@{}}
 \{1,2,3\}& \{2,3,5\}& \{0,3,9\}& \{1,3,4\}& \{2,6,7\}& \{0,3,4\}\\
 \{0,3,4\}& \{1,2,4\}& \{0,2,7\}& \{4,5,8\}& \{1,5,6\}& \{2,8,9\}\\
 \{2,8,9\}& \{4,7,9\}& \{1,6,8\}& \{2,5,9\}& \{4,7,8\}& \{1,5,7\}\\
 \{1,5,7\}& \{0,2,8\}& \{4,6,9\}& \{0,1,8\}& \{3,5,7\}& \{0,4,5\}\\
 \{0,4,5\}& \{1,7,9\}& \{3,6,9\}& \{5,8,9\}& \{0,1,9\}& \{3,7,8\}
\end{array}$
\end{description}}
\vbox{\begin{description}
\item[\T{10}{3}:] [(10)]\\

$\begin{array}{@{}*{5}{l@{\hspace{4pt}}}l@{}}
 \{1,2,3\}& \{2,3,5\}& \{0,3,9\}& \{1,3,4\}& \{0,2,6\}& \{0,3,4\}\\
 \{0,3,4\}& \{1,2,4\}& \{2,7,9\}& \{4,6,7\}& \{1,5,6\}& \{2,6,8\}\\
 \{2,6,8\}& \{4,7,8\}& \{1,6,8\}& \{2,8,9\}& \{4,8,9\}& \{1,5,7\}\\
 \{1,5,7\}& \{0,2,7\}& \{4,5,9\}& \{0,1,8\}& \{3,5,8\}& \{0,4,6\}\\
 \{0,4,6\}& \{1,7,9\}& \{3,6,9\}& \{0,5,7\}& \{0,1,9\}& \{3,6,7\}
\end{array}$
\end{description}}
\vbox{\begin{description}
\item[\T{10}{4}:] [(28)]\\

$\begin{array}{@{}*{5}{l@{\hspace{4pt}}}l@{}}
 \{1,2,3\}& \{2,3,5\}& \{3,4,9\}& \{1,3,4\}& \{2,6,7\}& \{0,3,6\}\\
 \{0,3,6\}& \{1,2,4\}& \{2,7,8\}& \{0,4,5\}& \{1,5,7\}& \{0,2,8\}\\
 \{0,2,8\}& \{4,6,7\}& \{1,5,6\}& \{2,5,9\}& \{4,5,8\}& \{1,7,9\}\\
 \{1,7,9\}& \{0,2,9\}& \{4,8,9\}& \{1,6,8\}& \{3,5,8\}& \{0,4,7\}\\
 \{0,4,7\}& \{0,1,9\}& \{3,6,9\}& \{0,5,6\}& \{0,1,8\}& \{0,3,7\}
\end{array}$
\end{description}}
\vbox{\begin{description}
\item[\T{10}{5}:] [(8)]\\

$\begin{array}{@{}*{5}{l@{\hspace{4pt}}}l@{}}
 \{1,2,3\}& \{2,3,5\}& \{3,7,9\}& \{1,3,4\}& \{2,6,8\}& \{0,3,4\}\\
 \{0,3,4\}& \{1,2,4\}& \{0,2,7\}& \{4,6,7\}& \{1,5,7\}& \{2,7,8\}\\
 \{2,7,8\}& \{4,7,8\}& \{1,5,6\}& \{2,6,9\}& \{4,8,9\}& \{1,7,9\}\\
 \{1,7,9\}& \{0,2,9\}& \{4,5,9\}& \{1,6,8\}& \{3,5,8\}& \{0,4,6\}\\
 \{0,4,6\}& \{0,1,9\}& \{3,6,9\}& \{0,5,6\}& \{0,1,8\}& \{3,6,7\}
\end{array}$
\end{description}}
\vbox{\begin{description}
\item[\T{10}{6}:] [(31)]\\

$\begin{array}{@{}*{5}{l@{\hspace{4pt}}}l@{}}
 \{1,2,3\}& \{1,3,4\}& \{1,2,4\}& \{1,5,6\}& \{1,6,8\}& \{1,5,7\}\\
 \{1,5,7\}& \{0,1,8\}& \{1,7,9\}& \{0,1,9\}& \{2,4,9\}& \{2,3,5\}\\
 \{2,3,5\}& \{2,5,6\}& \{2,6,7\}& \{2,7,8\}& \{0,2,9\}& \{0,2,8\}\\
 \{0,2,8\}& \{3,5,9\}& \{0,3,6\}& \{3,7,8\}& \{3,4,8\}& \{3,6,9\}\\
 \{3,6,9\}& \{0,3,7\}& \{0,4,5\}& \{4,6,7\}& \{4,7,9\}& \{4,5,8\}
\end{array}$
\end{description}}
\vbox{\begin{description}
\item[\T{10}{7}:] [(5)]\\

$\begin{array}{@{}*{5}{l@{\hspace{4pt}}}l@{}}
 \{1,2,3\}& \{2,3,5\}& \{3,4,9\}& \{1,3,4\}& \{2,6,8\}& \{0,3,7\}\\
 \{0,3,7\}& \{1,2,4\}& \{0,2,7\}& \{0,4,6\}& \{1,5,6\}& \{2,7,8\}\\
 \{2,7,8\}& \{4,7,8\}& \{1,6,7\}& \{2,6,9\}& \{4,6,8\}& \{1,5,7\}\\
 \{1,5,7\}& \{0,2,9\}& \{4,7,9\}& \{1,8,9\}& \{3,5,8\}& \{0,4,5\}\\
 \{0,4,5\}& \{0,1,9\}& \{3,6,9\}& \{5,7,9\}& \{0,1,8\}& \{3,6,7\}
\end{array}$
\end{description}}
\vbox{\begin{description}
\item[\T{10}{8}:] [(16)]\\

$\begin{array}{@{}*{5}{l@{\hspace{4pt}}}l@{}}
 \{1,2,3\}& \{2,3,5\}& \{3,4,9\}& \{1,3,4\}& \{2,6,9\}& \{0,3,8\}\\
 \{0,3,8\}& \{1,2,4\}& \{2,7,8\}& \{4,5,8\}& \{1,5,7\}& \{2,5,8\}\\
 \{2,5,8\}& \{0,4,7\}& \{1,5,6\}& \{0,2,9\}& \{4,7,8\}& \{1,6,7\}\\
 \{1,6,7\}& \{0,2,7\}& \{4,5,9\}& \{0,1,8\}& \{0,3,5\}& \{0,4,6\}\\
 \{0,4,6\}& \{1,8,9\}& \{3,6,7\}& \{0,5,6\}& \{0,1,9\}& \{3,7,9\}
\end{array}$
\end{description}}
\vbox{\begin{description}
\item[\T{10}{9}:] [(29)]\\

$\begin{array}{@{}*{5}{l@{\hspace{4pt}}}l@{}}
 \{1,2,3\}& \{1,3,4\}& \{1,2,4\}& \{1,5,6\}& \{1,6,8\}& \{1,5,7\}\\
 \{1,5,7\}& \{0,1,8\}& \{1,7,9\}& \{0,1,9\}& \{2,4,6\}& \{2,3,5\}\\
 \{2,3,5\}& \{0,2,6\}& \{2,7,9\}& \{2,7,8\}& \{2,8,9\}& \{0,2,5\}\\
 \{0,2,5\}& \{3,5,8\}& \{3,6,9\}& \{3,4,7\}& \{3,6,8\}& \{0,3,9\}\\
 \{0,3,9\}& \{0,3,7\}& \{0,4,5\}& \{4,6,7\}& \{4,8,9\}& \{4,5,9\}
\end{array}$
\end{description}}
\vbox{\begin{description}
\item[\T{10}{10}:] [(4)]\\

$\begin{array}{@{}*{5}{l@{\hspace{4pt}}}l@{}}
 \{1,2,3\}& \{2,3,5\}& \{3,6,9\}& \{1,3,4\}& \{2,6,8\}& \{0,3,4\}\\
 \{0,3,4\}& \{1,2,4\}& \{2,6,7\}& \{4,6,9\}& \{1,5,6\}& \{0,2,8\}\\
 \{0,2,8\}& \{4,7,8\}& \{1,6,7\}& \{2,7,9\}& \{4,6,8\}& \{1,5,7\}\\
 \{1,5,7\}& \{0,2,9\}& \{4,5,9\}& \{1,8,9\}& \{3,5,8\}& \{0,4,7\}\\
 \{0,4,7\}& \{0,1,9\}& \{0,3,6\}& \{0,5,7\}& \{0,1,8\}& \{3,7,9\}
\end{array}$
\end{description}}
\vbox{\begin{description}
\item[\T{10}{11}:] [(21)]\\

$\begin{array}{@{}*{5}{l@{\hspace{4pt}}}l@{}}
 \{1,2,3\}& \{2,3,5\}& \{3,6,9\}& \{1,3,4\}& \{2,6,9\}& \{0,3,4\}\\
 \{0,3,4\}& \{1,2,4\}& \{2,5,7\}& \{4,5,9\}& \{1,5,7\}& \{0,2,8\}\\
 \{0,2,8\}& \{4,7,8\}& \{1,5,6\}& \{2,8,9\}& \{4,5,8\}& \{1,7,9\}\\
 \{1,7,9\}& \{0,2,7\}& \{4,7,9\}& \{1,6,8\}& \{3,5,8\}& \{0,4,6\}\\
 \{0,4,6\}& \{0,1,9\}& \{3,6,7\}& \{0,5,6\}& \{0,1,8\}& \{0,3,7\}
\end{array}$
\end{description}}
\vbox{\begin{description}
\item[\T{10}{12}:] [(12)]\\

$\begin{array}{@{}*{5}{l@{\hspace{4pt}}}l@{}}
 \{1,2,3\}& \{2,3,5\}& \{3,8,9\}& \{1,3,4\}& \{2,6,9\}& \{0,3,6\}\\
 \{0,3,6\}& \{1,2,4\}& \{2,7,8\}& \{0,4,6\}& \{1,5,6\}& \{2,6,8\}\\
 \{2,6,8\}& \{4,7,9\}& \{1,6,8\}& \{0,2,9\}& \{4,7,8\}& \{1,5,7\}\\
 \{1,5,7\}& \{0,2,7\}& \{4,6,9\}& \{0,1,8\}& \{3,5,9\}& \{0,4,5\}\\
 \{0,4,5\}& \{1,7,9\}& \{3,6,7\}& \{5,6,7\}& \{0,1,9\}& \{0,3,7\}
\end{array}$
\end{description}}
\vbox{\begin{description}
\item[\T{10}{13}:] [(20)]\\

$\begin{array}{@{}*{5}{l@{\hspace{4pt}}}l@{}}
 \{1,2,3\}& \{2,3,5\}& \{3,6,9\}& \{1,3,4\}& \{2,6,9\}& \{0,3,9\}\\
 \{0,3,9\}& \{1,2,4\}& \{2,5,7\}& \{4,5,9\}& \{1,5,7\}& \{0,2,8\}\\
 \{0,2,8\}& \{0,4,7\}& \{1,5,6\}& \{2,8,9\}& \{4,5,8\}& \{1,7,9\}\\
 \{1,7,9\}& \{0,2,7\}& \{4,7,9\}& \{1,6,8\}& \{0,3,5\}& \{0,4,6\}\\
 \{0,4,6\}& \{0,1,9\}& \{3,6,7\}& \{0,5,6\}& \{0,1,8\}& \{3,7,8\}
\end{array}$
\end{description}}
\vbox{\begin{description}
\item[\T{10}{14}:] [(9)]\\

$\begin{array}{@{}*{5}{l@{\hspace{4pt}}}l@{}}
 \{1,2,3\}& \{2,3,5\}& \{3,6,9\}& \{1,3,4\}& \{2,6,8\}& \{0,3,9\}\\
 \{0,3,9\}& \{1,2,4\}& \{0,2,7\}& \{4,6,9\}& \{1,5,7\}& \{2,8,9\}\\
 \{2,8,9\}& \{4,6,7\}& \{1,5,6\}& \{2,7,9\}& \{0,4,8\}& \{1,7,9\}\\
 \{1,7,9\}& \{0,2,6\}& \{4,5,9\}& \{1,6,8\}& \{0,3,5\}& \{0,4,7\}\\
 \{0,4,7\}& \{0,1,9\}& \{3,6,7\}& \{0,5,6\}& \{0,1,8\}& \{3,7,8\}
\end{array}$
\end{description}}
\vbox{\begin{description}
\item[\T{10}{15}:] [(14)]\\

$\begin{array}{@{}*{5}{l@{\hspace{4pt}}}l@{}}
 \{1,2,3\}& \{2,3,5\}& \{3,6,9\}& \{1,3,4\}& \{2,6,8\}& \{0,3,8\}\\
 \{0,3,8\}& \{1,2,4\}& \{2,7,9\}& \{0,4,5\}& \{1,5,6\}& \{2,5,8\}\\
 \{2,5,8\}& \{4,7,8\}& \{1,6,7\}& \{0,2,9\}& \{4,8,9\}& \{1,5,7\}\\
 \{1,5,7\}& \{0,2,7\}& \{4,5,9\}& \{0,1,8\}& \{3,5,9\}& \{0,4,6\}\\
 \{0,4,6\}& \{1,8,9\}& \{0,3,6\}& \{0,5,7\}& \{0,1,9\}& \{3,4,7\}
\end{array}$
\end{description}}
\vbox{\begin{description}
\item[\T{10}{16}:] [(24)]\\

$\begin{array}{@{}*{5}{l@{\hspace{4pt}}}l@{}}
 \{1,2,3\}& \{2,3,5\}& \{3,8,9\}& \{1,3,4\}& \{2,6,7\}& \{0,3,6\}\\
 \{0,3,6\}& \{1,2,4\}& \{2,7,8\}& \{4,5,8\}& \{1,5,7\}& \{0,2,8\}\\
 \{0,2,8\}& \{0,4,7\}& \{1,5,6\}& \{2,5,9\}& \{4,8,9\}& \{1,7,9\}\\
 \{1,7,9\}& \{0,2,9\}& \{4,6,9\}& \{1,6,8\}& \{0,3,5\}& \{0,4,5\}\\
 \{0,4,5\}& \{0,1,9\}& \{3,6,9\}& \{5,6,8\}& \{0,1,8\}& \{3,4,7\}
\end{array}$
\end{description}}
\vbox{\begin{description}
\item[\T{10}{17}:] [(7)]\\

$\begin{array}{@{}*{5}{l@{\hspace{4pt}}}l@{}}
 \{1,2,3\}& \{2,3,5\}& \{3,4,9\}& \{1,3,4\}& \{2,6,9\}& \{0,3,6\}\\
 \{0,3,6\}& \{1,2,4\}& \{2,7,8\}& \{0,4,6\}& \{1,5,6\}& \{2,6,8\}\\
 \{2,6,8\}& \{4,6,7\}& \{1,6,8\}& \{0,2,9\}& \{4,7,8\}& \{1,5,7\}\\
 \{1,5,7\}& \{0,2,7\}& \{4,8,9\}& \{0,1,8\}& \{3,5,8\}& \{0,4,5\}\\
 \{0,4,5\}& \{1,7,9\}& \{3,6,7\}& \{5,7,9\}& \{0,1,9\}& \{0,3,7\}
\end{array}$
\end{description}}
\vbox{\begin{description}
\item[\T{10}{18}:] [(22)]\\

$\begin{array}{@{}*{5}{l@{\hspace{4pt}}}l@{}}
 \{1,2,3\}& \{2,3,5\}& \{3,6,9\}& \{1,3,4\}& \{2,6,7\}& \{0,3,7\}\\
 \{0,3,7\}& \{1,2,4\}& \{0,2,7\}& \{0,4,5\}& \{1,5,7\}& \{2,5,8\}\\
 \{2,5,8\}& \{4,7,8\}& \{1,5,6\}& \{2,8,9\}& \{4,8,9\}& \{1,7,9\}\\
 \{1,7,9\}& \{0,2,9\}& \{4,5,9\}& \{1,6,8\}& \{3,5,9\}& \{0,4,6\}\\
 \{0,4,6\}& \{0,1,9\}& \{3,6,8\}& \{0,5,6\}& \{0,1,8\}& \{3,4,7\}
\end{array}$
\end{description}}
\vbox{\begin{description}
\item[\T{10}{19}:] [(25)]\\

$\begin{array}{@{}*{5}{l@{\hspace{4pt}}}l@{}}
 \{1,2,3\}& \{2,3,5\}& \{3,8,9\}& \{1,3,4\}& \{2,6,7\}& \{0,3,4\}\\
 \{0,3,4\}& \{1,2,4\}& \{2,7,8\}& \{4,5,8\}& \{1,5,6\}& \{2,8,9\}\\
 \{2,8,9\}& \{4,7,9\}& \{1,6,8\}& \{0,2,9\}& \{4,7,8\}& \{1,5,7\}\\
 \{1,5,7\}& \{0,2,5\}& \{4,5,9\}& \{0,1,8\}& \{3,5,7\}& \{0,4,6\}\\
 \{0,4,6\}& \{1,7,9\}& \{3,6,9\}& \{5,6,9\}& \{0,1,9\}& \{0,3,7\}
\end{array}$
\end{description}}
\vbox{\begin{description}
\item[\T{10}{20}:] [(27)]\\

$\begin{array}{@{}*{5}{l@{\hspace{4pt}}}l@{}}
 \{1,2,3\}& \{2,3,5\}& \{3,7,9\}& \{1,3,4\}& \{2,6,7\}& \{0,3,6\}\\
 \{0,3,6\}& \{1,2,4\}& \{2,7,8\}& \{0,4,5\}& \{1,5,7\}& \{0,2,8\}\\
 \{0,2,8\}& \{4,5,7\}& \{1,5,6\}& \{2,5,9\}& \{4,8,9\}& \{1,7,9\}\\
 \{1,7,9\}& \{0,2,9\}& \{4,6,9\}& \{1,6,8\}& \{3,5,8\}& \{0,4,7\}\\
 \{0,4,7\}& \{0,1,9\}& \{3,6,9\}& \{0,5,6\}& \{0,1,8\}& \{0,3,7\}
\end{array}$
\end{description}}
\vbox{\begin{description}
\item[\T{10}{21}:] [(11)]\\

$\begin{array}{@{}*{5}{l@{\hspace{4pt}}}l@{}}
 \{1,2,3\}& \{2,3,5\}& \{0,3,9\}& \{1,3,4\}& \{0,2,6\}& \{0,3,4\}\\
 \{0,3,4\}& \{1,2,4\}& \{2,7,8\}& \{4,6,9\}& \{1,5,6\}& \{2,8,9\}\\
 \{2,8,9\}& \{4,6,7\}& \{1,6,8\}& \{2,6,9\}& \{4,7,8\}& \{1,5,7\}\\
 \{1,5,7\}& \{0,2,7\}& \{4,5,9\}& \{0,1,8\}& \{3,5,8\}& \{0,4,8\}\\
 \{0,4,8\}& \{1,7,9\}& \{3,6,7\}& \{0,5,7\}& \{0,1,9\}& \{3,7,9\}
\end{array}$
\end{description}}
\vbox{\begin{description}
\item[\T{10}{22}:] [(23)]\\

$\begin{array}{@{}*{5}{l@{\hspace{4pt}}}l@{}}
 \{1,2,3\}& \{2,3,5\}& \{3,6,9\}& \{1,3,4\}& \{2,6,7\}& \{0,3,7\}\\
 \{0,3,7\}& \{1,2,4\}& \{2,7,8\}& \{4,5,9\}& \{1,5,7\}& \{0,2,8\}\\
 \{0,2,8\}& \{0,4,7\}& \{1,5,6\}& \{2,5,9\}& \{4,5,8\}& \{1,7,9\}\\
 \{1,7,9\}& \{0,2,9\}& \{4,8,9\}& \{1,6,8\}& \{0,3,5\}& \{0,4,6\}\\
 \{0,4,6\}& \{0,1,9\}& \{3,6,8\}& \{0,5,6\}& \{0,1,8\}& \{3,4,7\}
\end{array}$
\end{description}}
\vbox{\begin{description}
\item[\T{10}{23}:] [(15)]\\

$\begin{array}{@{}*{5}{l@{\hspace{4pt}}}l@{}}
 \{1,2,3\}& \{2,3,5\}& \{3,6,9\}& \{1,3,4\}& \{2,6,9\}& \{0,3,8\}\\
 \{0,3,8\}& \{1,2,4\}& \{2,7,8\}& \{4,5,9\}& \{1,5,6\}& \{2,5,8\}\\
 \{2,5,8\}& \{0,4,7\}& \{1,6,7\}& \{0,2,9\}& \{4,6,8\}& \{1,5,7\}\\
 \{1,5,7\}& \{0,2,7\}& \{4,7,9\}& \{0,1,8\}& \{3,5,7\}& \{0,4,5\}\\
 \{0,4,5\}& \{1,8,9\}& \{0,3,6\}& \{5,8,9\}& \{0,1,9\}& \{3,7,9\}
\end{array}$
\end{description}}
\vbox{\begin{description}
\item[\T{10}{24}:] [(17)]\\

$\begin{array}{@{}*{5}{l@{\hspace{4pt}}}l@{}}
 \{1,2,3\}& \{2,3,5\}& \{3,4,9\}& \{1,3,4\}& \{2,6,9\}& \{0,3,8\}\\
 \{0,3,8\}& \{1,2,4\}& \{2,7,8\}& \{4,5,8\}& \{1,5,6\}& \{2,5,8\}\\
 \{2,5,8\}& \{0,4,7\}& \{1,6,7\}& \{0,2,9\}& \{4,6,8\}& \{1,5,7\}\\
 \{1,5,7\}& \{0,2,7\}& \{4,7,9\}& \{0,1,8\}& \{3,5,9\}& \{0,4,5\}\\
 \{0,4,5\}& \{1,8,9\}& \{0,3,6\}& \{5,7,9\}& \{0,1,9\}& \{3,6,7\}
\end{array}$
\end{description}}
\vbox{\begin{description}
\item[\T{10}{25}:] [(19)]\\

$\begin{array}{@{}*{5}{l@{\hspace{4pt}}}l@{}}
 \{1,2,3\}& \{2,3,5\}& \{3,7,9\}& \{1,3,4\}& \{0,2,6\}& \{0,3,4\}\\
 \{0,3,4\}& \{1,2,4\}& \{2,5,7\}& \{4,5,9\}& \{1,5,7\}& \{2,7,8\}\\
 \{2,7,8\}& \{4,7,8\}& \{1,5,6\}& \{2,8,9\}& \{4,5,8\}& \{1,7,9\}\\
 \{1,7,9\}& \{0,2,9\}& \{4,6,9\}& \{1,6,8\}& \{0,3,5\}& \{0,4,7\}\\
 \{0,4,7\}& \{0,1,9\}& \{3,6,8\}& \{5,6,9\}& \{0,1,8\}& \{3,6,7\}
\end{array}$
\end{description}}
\vbox{\begin{description}
\item[\T{10}{26}:] [(6)]\\

$\begin{array}{@{}*{5}{l@{\hspace{4pt}}}l@{}}
 \{1,2,3\}& \{2,3,5\}& \{3,4,9\}& \{1,3,4\}& \{2,6,8\}& \{0,3,6\}\\
 \{0,3,6\}& \{1,2,4\}& \{0,2,7\}& \{0,4,6\}& \{1,5,7\}& \{2,7,8\}\\
 \{2,7,8\}& \{4,7,8\}& \{1,5,6\}& \{2,6,9\}& \{4,6,8\}& \{1,6,7\}\\
 \{1,6,7\}& \{0,2,9\}& \{4,7,9\}& \{1,8,9\}& \{3,5,8\}& \{0,4,5\}\\
 \{0,4,5\}& \{0,1,9\}& \{3,6,7\}& \{5,6,9\}& \{0,1,8\}& \{3,7,9\}
\end{array}$
\end{description}}
\vbox{\begin{description}
\item[\T{10}{27}:] [(2)]\\

$\begin{array}{@{}*{5}{l@{\hspace{4pt}}}l@{}}
 \{1,2,3\}& \{2,5,8\}& \{3,8,9\}& \{1,3,4\}& \{2,5,6\}& \{0,3,6\}\\
 \{0,3,6\}& \{1,2,4\}& \{2,7,9\}& \{4,5,9\}& \{1,5,6\}& \{0,2,8\}\\
 \{0,2,8\}& \{0,4,6\}& \{1,6,7\}& \{2,6,9\}& \{4,7,8\}& \{1,5,7\}\\
 \{1,5,7\}& \{0,2,7\}& \{4,6,8\}& \{1,8,9\}& \{3,5,7\}& \{4,7,9\}\\
 \{4,7,9\}& \{0,1,9\}& \{3,6,9\}& \{0,4,5\}& \{0,1,8\}& \{0,3,7\}
\end{array}$
\end{description}}
\vbox{\begin{description}
\item[\T{10}{28}:] [(18)]\\

$\begin{array}{@{}*{5}{l@{\hspace{4pt}}}l@{}}
 \{1,2,3\}& \{2,3,5\}& \{3,4,9\}& \{1,3,4\}& \{2,6,8\}& \{0,3,6\}\\
 \{0,3,6\}& \{1,2,4\}& \{2,5,7\}& \{4,5,8\}& \{1,5,7\}& \{0,2,8\}\\
 \{0,2,8\}& \{0,4,7\}& \{1,5,6\}& \{2,7,9\}& \{4,7,8\}& \{1,7,9\}\\
 \{1,7,9\}& \{0,2,9\}& \{4,6,9\}& \{1,6,8\}& \{0,3,5\}& \{0,4,5\}\\
 \{0,4,5\}& \{0,1,9\}& \{3,6,7\}& \{5,6,9\}& \{0,1,8\}& \{3,7,8\}
\end{array}$
\end{description}}
\vbox{\begin{description}
\item[\T{10}{29}:] [(3)]\\

$\begin{array}{@{}*{5}{l@{\hspace{4pt}}}l@{}}
 \{1,2,3\}& \{0,2,5\}& \{0,3,9\}& \{1,3,4\}& \{2,5,6\}& \{0,3,5\}\\
 \{0,3,5\}& \{1,2,4\}& \{2,7,9\}& \{4,5,8\}& \{1,5,6\}& \{2,7,8\}\\
 \{2,7,8\}& \{4,6,7\}& \{1,6,8\}& \{2,6,9\}& \{0,4,7\}& \{1,5,7\}\\
 \{1,5,7\}& \{0,2,8\}& \{4,8,9\}& \{0,1,8\}& \{3,5,7\}& \{4,5,9\}\\
 \{4,5,9\}& \{1,7,9\}& \{3,6,9\}& \{0,4,6\}& \{0,1,9\}& \{3,7,8\}
\end{array}$
\end{description}}
\vbox{\begin{description}
\item[\T{10}{30}:] [(13)]\\

$\begin{array}{@{}*{5}{l@{\hspace{4pt}}}l@{}}
 \{1,2,3\}& \{2,3,5\}& \{3,7,9\}& \{1,3,4\}& \{0,2,6\}& \{0,3,9\}\\
 \{0,3,9\}& \{1,2,4\}& \{2,7,8\}& \{4,6,9\}& \{1,5,7\}& \{2,8,9\}\\
 \{2,8,9\}& \{0,4,7\}& \{1,5,6\}& \{2,6,9\}& \{4,7,8\}& \{1,7,9\}\\
 \{1,7,9\}& \{0,2,7\}& \{4,5,9\}& \{1,6,8\}& \{0,3,5\}& \{0,4,6\}\\
 \{0,4,6\}& \{0,1,9\}& \{3,6,8\}& \{5,6,7\}& \{0,1,8\}& \{3,6,7\}
\end{array}$
\end{description}}
\vbox{\begin{description}
\item[\T{10}{31}:] [(30)]\\

$\begin{array}{@{}*{5}{l@{\hspace{4pt}}}l@{}}
 \{1,2,3\}& \{1,3,4\}& \{1,2,4\}& \{1,5,6\}& \{1,6,8\}& \{1,5,7\}\\
 \{1,5,7\}& \{0,1,8\}& \{1,7,9\}& \{0,1,9\}& \{2,4,6\}& \{2,3,5\}\\
 \{2,3,5\}& \{2,6,9\}& \{0,2,7\}& \{2,7,8\}& \{2,8,9\}& \{0,2,5\}\\
 \{0,2,5\}& \{3,5,8\}& \{3,6,7\}& \{3,7,9\}& \{3,6,8\}& \{0,3,9\}\\
 \{0,3,9\}& \{0,3,4\}& \{4,5,7\}& \{4,7,8\}& \{4,8,9\}& \{4,5,9\}
\end{array}$
\end{description}}
\vbox{\begin{description}
\item[\T{10}{32}:] [(32)]\\

$\begin{array}{@{}*{5}{l@{\hspace{4pt}}}l@{}}
 \{1,2,3\}& \{1,3,5\}& \{1,2,4\}& \{1,4,5\}& \{1,6,8\}& \{1,6,7\}\\
 \{1,6,7\}& \{0,1,8\}& \{1,7,9\}& \{0,1,9\}& \{2,4,7\}& \{2,5,6\}\\
 \{2,5,6\}& \{2,3,6\}& \{2,7,8\}& \{2,8,9\}& \{0,2,9\}& \{0,2,5\}\\
 \{0,2,5\}& \{3,4,9\}& \{0,3,6\}& \{3,5,7\}& \{3,7,8\}& \{3,8,9\}\\
 \{3,8,9\}& \{0,3,4\}& \{4,5,8\}& \{4,6,9\}& \{4,6,8\}& \{0,4,7\}
\end{array}$
\end{description}}

\vbox{\begin{description}
\item[\T{10}{33}:] [(1.41)]\\

$\begin{array}{@{}*{5}{l@{\hspace{4pt}}}l@{}}
 \{1,2,3\}& \{1,2,3\}& \{1,4,5\}& \{1,5,6\}& \{1,6,7\}& \{1,7,8\}\\
 \{1,7,8\}& \{1,8,9\}& \{0,1,9\}& \{0,1,4\}& \{2,4,6\}& \{2,5,9\}\\
 \{2,5,9\}& \{2,6,8\}& \{0,2,7\}& \{2,4,8\}& \{2,7,9\}& \{0,2,5\}\\
 \{0,2,5\}& \{3,4,7\}& \{3,4,5\}& \{3,6,9\}& \{3,5,7\}& \{0,3,8\}\\
 \{0,3,8\}& \{3,8,9\}& \{0,3,6\}& \{4,7,8\}& \{4,6,9\}& \{0,4,9\}
\end{array}$
\end{description}}
\vbox{\begin{description}
\item[\T{10}{34}:] [(1.4)]\\

$\begin{array}{@{}*{5}{l@{\hspace{4pt}}}l@{}}
 \{1,2,3\}& \{1,2,3\}& \{1,4,5\}& \{1,5,6\}& \{1,6,7\}& \{1,7,8\}\\
 \{1,7,8\}& \{1,8,9\}& \{0,1,9\}& \{0,1,4\}& \{2,4,6\}& \{2,5,8\}\\
 \{2,5,8\}& \{2,5,6\}& \{2,7,9\}& \{0,2,8\}& \{2,4,9\}& \{0,2,7\}\\
 \{0,2,7\}& \{3,4,9\}& \{0,3,5\}& \{3,6,8\}& \{3,4,7\}& \{3,7,8\}\\
 \{3,7,8\}& \{3,5,9\}& \{0,3,6\}& \{4,5,7\}& \{4,6,8\}& \{0,4,8\}
\end{array}$
\end{description}}
\vbox{\begin{description}
\item[\T{10}{35}:] [(1.9)]\\

$\begin{array}{@{}*{5}{l@{\hspace{4pt}}}l@{}}
 \{1,2,3\}& \{1,2,3\}& \{1,4,5\}& \{1,5,6\}& \{1,6,7\}& \{1,7,8\}\\
 \{1,7,8\}& \{1,8,9\}& \{0,1,9\}& \{0,1,4\}& \{2,4,6\}& \{0,2,5\}\\
 \{0,2,5\}& \{2,6,8\}& \{2,7,9\}& \{2,5,8\}& \{2,4,9\}& \{0,2,7\}\\
 \{0,2,7\}& \{3,4,7\}& \{3,5,8\}& \{3,6,9\}& \{0,3,7\}& \{3,4,8\}\\
 \{3,4,8\}& \{3,5,9\}& \{0,3,6\}& \{4,7,8\}& \{4,6,9\}& \{0,4,5\}
\end{array}$
\end{description}}
\vbox{\begin{description}
\item[\T{10}{36}:] [(1.36)]\\

$\begin{array}{@{}*{5}{l@{\hspace{4pt}}}l@{}}
 \{1,2,3\}& \{1,2,3\}& \{1,4,5\}& \{1,5,6\}& \{1,4,6\}& \{1,7,8\}\\
 \{1,7,8\}& \{1,8,9\}& \{0,1,9\}& \{0,1,7\}& \{2,4,7\}& \{0,2,5\}\\
 \{0,2,5\}& \{2,6,8\}& \{2,6,7\}& \{2,5,8\}& \{2,4,9\}& \{0,2,9\}\\
 \{0,2,9\}& \{3,4,8\}& \{3,5,7\}& \{3,6,9\}& \{3,4,7\}& \{0,3,8\}\\
 \{0,3,8\}& \{3,5,9\}& \{0,3,6\}& \{0,4,6\}& \{4,8,9\}& \{0,4,5\}
\end{array}$
\end{description}}
\vbox{\begin{description}
\item[\T{10}{37}:] [(1.45)]\\

$\begin{array}{@{}*{5}{l@{\hspace{4pt}}}l@{}}
 \{1,2,3\}& \{1,2,3\}& \{1,4,5\}& \{1,5,6\}& \{1,4,6\}& \{1,7,8\}\\
 \{1,7,8\}& \{1,8,9\}& \{0,1,9\}& \{0,1,7\}& \{2,4,7\}& \{2,5,8\}\\
 \{2,5,8\}& \{2,5,6\}& \{2,7,9\}& \{0,2,8\}& \{2,4,9\}& \{0,2,6\}\\
 \{0,2,6\}& \{3,4,9\}& \{0,3,5\}& \{3,6,8\}& \{3,6,7\}& \{3,7,8\}\\
 \{3,7,8\}& \{3,5,9\}& \{0,3,4\}& \{4,6,7\}& \{4,5,8\}& \{0,4,8\}
\end{array}$
\end{description}}
\vbox{\begin{description}
\item[\T{10}{38}:] [(1.15)]\\

$\begin{array}{@{}*{5}{l@{\hspace{4pt}}}l@{}}
 \{1,2,3\}& \{1,2,3\}& \{1,4,5\}& \{1,5,6\}& \{1,6,7\}& \{1,7,8\}\\
 \{1,7,8\}& \{1,8,9\}& \{0,1,9\}& \{0,1,4\}& \{2,4,6\}& \{0,2,5\}\\
 \{0,2,5\}& \{2,6,9\}& \{2,4,7\}& \{2,7,8\}& \{2,5,9\}& \{0,2,8\}\\
 \{0,2,8\}& \{3,4,8\}& \{3,5,7\}& \{3,4,6\}& \{0,3,7\}& \{3,5,8\}\\
 \{3,5,8\}& \{3,6,9\}& \{0,3,9\}& \{4,7,9\}& \{4,8,9\}& \{0,4,5\}
\end{array}$
\end{description}}
\vbox{\begin{description}
\item[\T{10}{39}:] [(1.16)]\\

$\begin{array}{@{}*{5}{l@{\hspace{4pt}}}l@{}}
 \{1,2,3\}& \{1,2,3\}& \{1,4,5\}& \{1,5,6\}& \{1,4,6\}& \{1,7,8\}\\
 \{1,7,8\}& \{1,8,9\}& \{0,1,9\}& \{0,1,7\}& \{2,4,6\}& \{2,5,8\}\\
 \{2,5,8\}& \{2,6,7\}& \{2,5,7\}& \{0,2,8\}& \{2,4,9\}& \{0,2,9\}\\
 \{0,2,9\}& \{3,4,8\}& \{0,3,5\}& \{3,6,9\}& \{3,4,7\}& \{3,5,8\}\\
 \{3,5,8\}& \{3,7,9\}& \{0,3,6\}& \{0,4,7\}& \{4,8,9\}& \{0,4,5\}
\end{array}$
\end{description}}
\vbox{\begin{description}
\item[\T{10}{40}:] [(1.40)]\\

$\begin{array}{@{}*{5}{l@{\hspace{4pt}}}l@{}}
 \{1,2,3\}& \{1,2,3\}& \{1,4,5\}& \{1,5,6\}& \{1,6,7\}& \{1,7,8\}\\
 \{1,7,8\}& \{1,8,9\}& \{0,1,9\}& \{0,1,4\}& \{2,4,6\}& \{2,5,7\}\\
 \{2,5,7\}& \{2,6,8\}& \{0,2,7\}& \{2,4,8\}& \{2,5,9\}& \{0,2,9\}\\
 \{0,2,9\}& \{3,4,8\}& \{0,3,5\}& \{3,5,6\}& \{3,7,9\}& \{3,7,8\}\\
 \{3,7,8\}& \{3,4,9\}& \{0,3,6\}& \{4,5,7\}& \{4,6,9\}& \{0,4,7\}
\end{array}$
\end{description}}
\vbox{\begin{description}
\item[\T{10}{41}:] [(1.32)]\\

$\begin{array}{@{}*{5}{l@{\hspace{4pt}}}l@{}}
 \{1,2,3\}& \{1,2,3\}& \{1,4,5\}& \{1,5,6\}& \{1,4,6\}& \{1,7,8\}\\
 \{1,7,8\}& \{1,8,9\}& \{0,1,9\}& \{0,1,7\}& \{2,4,7\}& \{2,5,8\}\\
 \{2,5,8\}& \{0,2,6\}& \{2,5,7\}& \{2,6,8\}& \{2,4,9\}& \{0,2,9\}\\
 \{0,2,9\}& \{0,3,4\}& \{3,5,7\}& \{3,6,8\}& \{3,7,9\}& \{3,4,8\}\\
 \{3,4,8\}& \{3,6,9\}& \{0,3,5\}& \{4,6,7\}& \{0,4,8\}& \{4,5,9\}
\end{array}$
\end{description}}
\vbox{\begin{description}
\item[\T{10}{42}:] [(1.17)]\\

$\begin{array}{@{}*{5}{l@{\hspace{4pt}}}l@{}}
 \{1,2,3\}& \{1,2,3\}& \{1,4,5\}& \{1,5,6\}& \{1,4,6\}& \{1,7,8\}\\
 \{1,7,8\}& \{1,8,9\}& \{0,1,9\}& \{0,1,7\}& \{2,4,6\}& \{2,5,8\}\\
 \{2,5,8\}& \{2,6,7\}& \{2,5,7\}& \{0,2,8\}& \{2,4,9\}& \{0,2,9\}\\
 \{0,2,9\}& \{0,3,4\}& \{3,5,7\}& \{3,6,8\}& \{3,7,9\}& \{3,4,8\}\\
 \{3,4,8\}& \{3,6,9\}& \{0,3,5\}& \{0,4,7\}& \{4,7,8\}& \{4,5,9\}
\end{array}$
\end{description}}
\vbox{\begin{description}
\item[\T{10}{43}:] [(1.5)]\\

$\begin{array}{@{}*{5}{l@{\hspace{4pt}}}l@{}}
 \{1,2,3\}& \{1,2,3\}& \{1,4,5\}& \{1,5,6\}& \{1,6,7\}& \{1,7,8\}\\
 \{1,7,8\}& \{1,8,9\}& \{0,1,9\}& \{0,1,4\}& \{2,4,6\}& \{2,5,8\}\\
 \{2,5,8\}& \{2,5,6\}& \{2,7,9\}& \{0,2,8\}& \{2,4,9\}& \{0,2,7\}\\
 \{0,2,7\}& \{3,4,9\}& \{0,3,5\}& \{3,6,8\}& \{3,4,7\}& \{3,7,8\}\\
 \{3,7,8\}& \{3,5,9\}& \{0,3,6\}& \{4,6,7\}& \{4,5,8\}& \{0,4,8\}
\end{array}$
\end{description}}
\vbox{\begin{description}
\item[\T{10}{44}:] [(1.18)]\\

$\begin{array}{@{}*{5}{l@{\hspace{4pt}}}l@{}}
 \{1,2,3\}& \{1,2,3\}& \{1,4,5\}& \{1,5,6\}& \{1,4,6\}& \{1,7,8\}\\
 \{1,7,8\}& \{1,8,9\}& \{0,1,9\}& \{0,1,7\}& \{2,4,6\}& \{2,5,8\}\\
 \{2,5,8\}& \{2,6,7\}& \{2,5,7\}& \{0,2,8\}& \{2,4,9\}& \{0,2,9\}\\
 \{0,2,9\}& \{0,3,4\}& \{3,5,9\}& \{3,6,8\}& \{3,4,7\}& \{3,7,8\}\\
 \{3,7,8\}& \{3,6,9\}& \{0,3,5\}& \{0,4,7\}& \{4,5,8\}& \{4,8,9\}
\end{array}$
\end{description}}
\vbox{\begin{description}
\item[\T{10}{45}:] [(1.22)]\\

$\begin{array}{@{}*{5}{l@{\hspace{4pt}}}l@{}}
 \{1,2,3\}& \{1,2,3\}& \{1,4,5\}& \{1,5,6\}& \{1,4,6\}& \{1,7,8\}\\
 \{1,7,8\}& \{1,8,9\}& \{0,1,9\}& \{0,1,7\}& \{2,4,6\}& \{2,5,9\}\\
 \{2,5,9\}& \{2,6,7\}& \{2,5,7\}& \{0,2,8\}& \{2,8,9\}& \{0,2,4\}\\
 \{0,2,4\}& \{3,4,9\}& \{3,5,7\}& \{3,6,8\}& \{0,3,7\}& \{3,4,8\}\\
 \{3,4,8\}& \{3,5,9\}& \{0,3,6\}& \{4,7,9\}& \{4,7,8\}& \{0,4,5\}
\end{array}$
\end{description}}
\vbox{\begin{description}
\item[\T{10}{46}:] [(1.30)]\\

$\begin{array}{@{}*{5}{l@{\hspace{4pt}}}l@{}}
 \{1,2,3\}& \{1,2,3\}& \{1,4,5\}& \{1,5,6\}& \{1,4,6\}& \{1,7,8\}\\
 \{1,7,8\}& \{1,8,9\}& \{0,1,9\}& \{0,1,7\}& \{2,4,7\}& \{2,5,8\}\\
 \{2,5,8\}& \{0,2,6\}& \{2,5,7\}& \{2,6,8\}& \{2,4,9\}& \{0,2,9\}\\
 \{0,2,9\}& \{3,4,8\}& \{0,3,5\}& \{3,6,9\}& \{3,6,7\}& \{3,7,8\}\\
 \{3,7,8\}& \{3,5,9\}& \{0,3,4\}& \{4,5,6\}& \{4,7,9\}& \{0,4,8\}
\end{array}$
\end{description}}
\vbox{\begin{description}
\item[\T{10}{47}:] [(1.33)]\\

$\begin{array}{@{}*{5}{l@{\hspace{4pt}}}l@{}}
 \{1,2,3\}& \{1,2,3\}& \{1,4,5\}& \{1,5,6\}& \{1,4,6\}& \{1,7,8\}\\
 \{1,7,8\}& \{1,8,9\}& \{0,1,9\}& \{0,1,7\}& \{2,4,7\}& \{2,5,8\}\\
 \{2,5,8\}& \{0,2,6\}& \{2,5,7\}& \{2,6,8\}& \{2,4,9\}& \{0,2,9\}\\
 \{0,2,9\}& \{0,3,4\}& \{3,5,9\}& \{3,6,8\}& \{3,4,7\}& \{3,7,8\}\\
 \{3,7,8\}& \{3,6,9\}& \{0,3,5\}& \{4,5,6\}& \{0,4,8\}& \{4,8,9\}
\end{array}$
\end{description}}
\vbox{\begin{description}
\item[\T{10}{48}:] [(1.21)]\\

$\begin{array}{@{}*{5}{l@{\hspace{4pt}}}l@{}}
 \{1,2,3\}& \{1,2,3\}& \{1,4,5\}& \{1,5,6\}& \{1,4,6\}& \{1,7,8\}\\
 \{1,7,8\}& \{1,8,9\}& \{0,1,9\}& \{0,1,7\}& \{2,4,6\}& \{2,5,9\}\\
 \{2,5,9\}& \{2,6,7\}& \{2,5,7\}& \{0,2,8\}& \{2,8,9\}& \{0,2,4\}\\
 \{0,2,4\}& \{3,4,9\}& \{0,3,5\}& \{3,6,8\}& \{3,4,7\}& \{3,7,8\}\\
 \{3,7,8\}& \{3,5,9\}& \{0,3,6\}& \{4,7,9\}& \{4,5,8\}& \{0,4,8\}
\end{array}$
\end{description}}
\vbox{\begin{description}
\item[\T{10}{49}:] [(1.13)]\\

$\begin{array}{@{}*{5}{l@{\hspace{4pt}}}l@{}}
 \{1,2,3\}& \{1,2,3\}& \{1,4,5\}& \{1,5,6\}& \{1,6,7\}& \{1,7,8\}\\
 \{1,7,8\}& \{1,8,9\}& \{0,1,9\}& \{0,1,4\}& \{2,4,6\}& \{2,5,8\}\\
 \{2,5,8\}& \{2,6,9\}& \{2,4,7\}& \{0,2,8\}& \{2,5,9\}& \{0,2,7\}\\
 \{0,2,7\}& \{3,4,7\}& \{3,4,5\}& \{0,3,6\}& \{3,7,9\}& \{3,5,8\}\\
 \{3,5,8\}& \{3,6,9\}& \{0,3,8\}& \{4,6,8\}& \{4,8,9\}& \{0,4,9\}
\end{array}$
\end{description}}
\vbox{\begin{description}
\item[\T{10}{50}:] [(1.39)]\\

$\begin{array}{@{}*{5}{l@{\hspace{4pt}}}l@{}}
 \{1,2,3\}& \{1,2,3\}& \{1,4,5\}& \{1,5,6\}& \{1,6,7\}& \{1,7,8\}\\
 \{1,7,8\}& \{1,8,9\}& \{0,1,9\}& \{0,1,4\}& \{2,4,6\}& \{2,5,7\}\\
 \{2,5,7\}& \{2,6,8\}& \{0,2,7\}& \{2,4,8\}& \{2,5,9\}& \{0,2,9\}\\
 \{0,2,9\}& \{3,4,9\}& \{3,5,8\}& \{3,4,6\}& \{3,5,7\}& \{0,3,8\}\\
 \{0,3,8\}& \{3,6,9\}& \{0,3,7\}& \{4,7,9\}& \{4,7,8\}& \{0,4,5\}
\end{array}$
\end{description}}
\vbox{\begin{description}
\item[\T{10}{51}:] [(1.19)]\\

$\begin{array}{@{}*{5}{l@{\hspace{4pt}}}l@{}}
 \{1,2,3\}& \{1,2,3\}& \{1,4,5\}& \{1,5,6\}& \{1,4,6\}& \{1,7,8\}\\
 \{1,7,8\}& \{1,8,9\}& \{0,1,9\}& \{0,1,7\}& \{2,4,6\}& \{2,5,9\}\\
 \{2,5,9\}& \{2,6,7\}& \{2,5,7\}& \{0,2,8\}& \{2,8,9\}& \{0,2,4\}\\
 \{0,2,4\}& \{3,4,7\}& \{3,5,8\}& \{3,6,9\}& \{0,3,7\}& \{3,6,8\}\\
 \{3,6,8\}& \{3,4,9\}& \{0,3,5\}& \{4,5,8\}& \{4,7,9\}& \{0,4,8\}
\end{array}$
\end{description}}
\vbox{\begin{description}
\item[\T{10}{52}:] [(1.10)]\\

$\begin{array}{@{}*{5}{l@{\hspace{4pt}}}l@{}}
 \{1,2,3\}& \{1,2,3\}& \{1,4,5\}& \{1,5,6\}& \{1,6,7\}& \{1,7,8\}\\
 \{1,7,8\}& \{1,8,9\}& \{0,1,9\}& \{0,1,4\}& \{2,4,6\}& \{0,2,5\}\\
 \{0,2,5\}& \{2,6,8\}& \{2,7,9\}& \{2,5,8\}& \{2,4,9\}& \{0,2,7\}\\
 \{0,2,7\}& \{3,4,7\}& \{3,4,5\}& \{3,6,9\}& \{0,3,7\}& \{3,6,8\}\\
 \{3,6,8\}& \{3,5,9\}& \{0,3,8\}& \{4,7,8\}& \{4,8,9\}& \{0,4,6\}
\end{array}$
\end{description}}
\vbox{\begin{description}
\item[\T{10}{53}:] [(1.31)]\\

$\begin{array}{@{}*{5}{l@{\hspace{4pt}}}l@{}}
 \{1,2,3\}& \{1,2,3\}& \{1,4,5\}& \{1,5,6\}& \{1,4,6\}& \{1,7,8\}\\
 \{1,7,8\}& \{1,8,9\}& \{0,1,9\}& \{0,1,7\}& \{2,4,7\}& \{2,5,8\}\\
 \{2,5,8\}& \{0,2,6\}& \{2,5,7\}& \{2,6,8\}& \{2,4,9\}& \{0,2,9\}\\
 \{0,2,9\}& \{3,4,9\}& \{3,5,7\}& \{3,6,8\}& \{0,3,7\}& \{3,4,8\}\\
 \{3,4,8\}& \{3,5,9\}& \{0,3,6\}& \{4,6,7\}& \{0,4,8\}& \{0,4,5\}
\end{array}$
\end{description}}
\vbox{\begin{description}
\item[\T{10}{54}:] [(1.34)]\\

$\begin{array}{@{}*{5}{l@{\hspace{4pt}}}l@{}}
 \{1,2,3\}& \{1,2,3\}& \{1,4,5\}& \{1,5,6\}& \{1,4,6\}& \{1,7,8\}\\
 \{1,7,8\}& \{1,8,9\}& \{0,1,9\}& \{0,1,7\}& \{2,4,7\}& \{2,5,8\}\\
 \{2,5,8\}& \{2,6,9\}& \{2,5,7\}& \{0,2,8\}& \{2,4,9\}& \{0,2,6\}\\
 \{0,2,6\}& \{3,4,8\}& \{0,3,5\}& \{3,6,7\}& \{3,7,9\}& \{3,6,8\}\\
 \{3,6,8\}& \{3,5,9\}& \{0,3,4\}& \{4,5,6\}& \{4,7,9\}& \{0,4,8\}
\end{array}$
\end{description}}
\vbox{\begin{description}
\item[\T{10}{55}:] [(1.27)]\\

$\begin{array}{@{}*{5}{l@{\hspace{4pt}}}l@{}}
 \{1,2,3\}& \{1,2,3\}& \{1,4,5\}& \{1,5,6\}& \{1,4,6\}& \{1,7,8\}\\
 \{1,7,8\}& \{1,8,9\}& \{0,1,9\}& \{0,1,7\}& \{2,4,6\}& \{0,2,5\}\\
 \{0,2,5\}& \{2,6,7\}& \{2,7,9\}& \{2,4,8\}& \{2,5,9\}& \{0,2,8\}\\
 \{0,2,8\}& \{3,4,8\}& \{3,5,7\}& \{0,3,6\}& \{3,6,7\}& \{3,5,8\}\\
 \{3,5,8\}& \{3,4,9\}& \{0,3,9\}& \{0,4,7\}& \{4,7,9\}& \{0,4,5\}
\end{array}$
\end{description}}
\vbox{\begin{description}
\item[\T{10}{56}:] [(1.7)]\\

$\begin{array}{@{}*{5}{l@{\hspace{4pt}}}l@{}}
 \{1,2,3\}& \{1,2,3\}& \{1,4,5\}& \{1,5,6\}& \{1,6,7\}& \{1,7,8\}\\
 \{1,7,8\}& \{1,8,9\}& \{0,1,9\}& \{0,1,4\}& \{2,4,6\}& \{2,5,9\}\\
 \{2,5,9\}& \{2,5,6\}& \{2,4,7\}& \{0,2,8\}& \{2,8,9\}& \{0,2,7\}\\
 \{0,2,7\}& \{0,3,4\}& \{3,5,8\}& \{3,6,9\}& \{3,5,7\}& \{3,4,8\}\\
 \{3,4,8\}& \{3,7,9\}& \{0,3,6\}& \{4,7,9\}& \{4,5,8\}& \{4,6,9\}
\end{array}$
\end{description}}
\vbox{\begin{description}
\item[\T{10}{57}:] [(1.23)]\\

$\begin{array}{@{}*{5}{l@{\hspace{4pt}}}l@{}}
 \{1,2,3\}& \{1,2,3\}& \{1,4,5\}& \{1,5,6\}& \{1,4,6\}& \{1,7,8\}\\
 \{1,7,8\}& \{1,8,9\}& \{0,1,9\}& \{0,1,7\}& \{2,4,6\}& \{0,2,5\}\\
 \{0,2,5\}& \{2,6,7\}& \{2,5,7\}& \{2,4,8\}& \{2,8,9\}& \{0,2,9\}\\
 \{0,2,9\}& \{0,3,4\}& \{3,5,7\}& \{3,6,8\}& \{3,7,9\}& \{3,5,8\}\\
 \{3,5,8\}& \{3,4,9\}& \{0,3,6\}& \{0,4,7\}& \{4,7,8\}& \{4,5,9\}
\end{array}$
\end{description}}
\vbox{\begin{description}
\item[\T{10}{58}:] [(1.38)]\\

$\begin{array}{@{}*{5}{l@{\hspace{4pt}}}l@{}}
 \{1,2,3\}& \{1,2,3\}& \{1,4,5\}& \{1,5,6\}& \{1,6,7\}& \{1,7,8\}\\
 \{1,7,8\}& \{1,8,9\}& \{0,1,9\}& \{0,1,4\}& \{2,4,6\}& \{2,5,7\}\\
 \{2,5,7\}& \{2,6,8\}& \{0,2,7\}& \{2,4,8\}& \{2,5,9\}& \{0,2,9\}\\
 \{0,2,9\}& \{3,4,7\}& \{3,5,8\}& \{3,5,6\}& \{3,7,9\}& \{0,3,8\}\\
 \{0,3,8\}& \{3,4,9\}& \{0,3,6\}& \{4,7,8\}& \{4,6,9\}& \{0,4,5\}
\end{array}$
\end{description}}
\vbox{\begin{description}
\item[\T{10}{59}:] [(1.29)]\\

$\begin{array}{@{}*{5}{l@{\hspace{4pt}}}l@{}}
 \{1,2,3\}& \{1,2,3\}& \{1,4,5\}& \{1,5,6\}& \{1,4,6\}& \{1,7,8\}\\
 \{1,7,8\}& \{1,8,9\}& \{0,1,9\}& \{0,1,7\}& \{2,4,6\}& \{2,5,9\}\\
 \{2,5,9\}& \{2,6,7\}& \{0,2,7\}& \{2,4,8\}& \{2,8,9\}& \{0,2,5\}\\
 \{0,2,5\}& \{0,3,4\}& \{3,5,7\}& \{3,6,8\}& \{3,4,7\}& \{3,5,8\}\\
 \{3,5,8\}& \{3,6,9\}& \{0,3,9\}& \{4,7,9\}& \{0,4,8\}& \{4,5,9\}
\end{array}$
\end{description}}
\vbox{\begin{description}
\item[\T{10}{60}:] [(1.2)]\\

$\begin{array}{@{}*{5}{l@{\hspace{4pt}}}l@{}}
 \{1,2,3\}& \{1,2,3\}& \{1,4,5\}& \{1,5,6\}& \{1,6,7\}& \{1,7,8\}\\
 \{1,7,8\}& \{1,8,9\}& \{0,1,9\}& \{0,1,4\}& \{2,4,6\}& \{2,5,8\}\\
 \{2,5,8\}& \{2,5,6\}& \{2,7,9\}& \{0,2,8\}& \{2,4,9\}& \{0,2,7\}\\
 \{0,2,7\}& \{3,4,7\}& \{0,3,5\}& \{3,6,8\}& \{3,5,7\}& \{3,4,8\}\\
 \{3,4,8\}& \{3,6,9\}& \{0,3,9\}& \{4,7,8\}& \{4,5,9\}& \{0,4,6\}
\end{array}$
\end{description}}
\vbox{\begin{description}
\item[\T{10}{61}:] [(1.6)]\\

$\begin{array}{@{}*{5}{l@{\hspace{4pt}}}l@{}}
 \{1,2,3\}& \{1,2,3\}& \{1,4,5\}& \{1,5,6\}& \{1,6,7\}& \{1,7,8\}\\
 \{1,7,8\}& \{1,8,9\}& \{0,1,9\}& \{0,1,4\}& \{2,4,6\}& \{2,5,9\}\\
 \{2,5,9\}& \{2,5,6\}& \{2,4,7\}& \{0,2,8\}& \{2,8,9\}& \{0,2,7\}\\
 \{0,2,7\}& \{3,4,7\}& \{0,3,5\}& \{3,6,8\}& \{3,7,9\}& \{3,4,8\}\\
 \{3,4,8\}& \{3,5,9\}& \{0,3,6\}& \{4,5,8\}& \{4,6,9\}& \{0,4,9\}
\end{array}$
\end{description}}
\vbox{\begin{description}
\item[\T{10}{62}:] [(1.1)]\\

$\begin{array}{@{}*{5}{l@{\hspace{4pt}}}l@{}}
 \{1,2,3\}& \{1,2,3\}& \{1,4,5\}& \{1,5,6\}& \{1,6,7\}& \{1,7,8\}\\
 \{1,7,8\}& \{1,8,9\}& \{0,1,9\}& \{0,1,4\}& \{2,4,6\}& \{2,5,8\}\\
 \{2,5,8\}& \{2,5,6\}& \{2,7,9\}& \{0,2,8\}& \{2,4,9\}& \{0,2,7\}\\
 \{0,2,7\}& \{3,4,7\}& \{3,5,9\}& \{0,3,6\}& \{3,5,7\}& \{3,4,8\}\\
 \{3,4,8\}& \{3,6,9\}& \{0,3,8\}& \{4,5,8\}& \{4,7,9\}& \{0,4,6\}
\end{array}$
\end{description}}
\vbox{\begin{description}
\item[\T{10}{63}:] [(1.24)]\\

$\begin{array}{@{}*{5}{l@{\hspace{4pt}}}l@{}}
 \{1,2,3\}& \{1,2,3\}& \{1,4,5\}& \{1,5,6\}& \{1,4,6\}& \{1,7,8\}\\
 \{1,7,8\}& \{1,8,9\}& \{0,1,9\}& \{0,1,7\}& \{2,4,6\}& \{2,5,8\}\\
 \{2,5,8\}& \{2,6,7\}& \{2,7,9\}& \{0,2,8\}& \{2,5,9\}& \{0,2,4\}\\
 \{0,2,4\}& \{3,4,7\}& \{3,5,9\}& \{0,3,6\}& \{3,5,7\}& \{3,4,8\}\\
 \{3,4,8\}& \{3,6,9\}& \{0,3,8\}& \{4,5,8\}& \{4,7,9\}& \{0,4,9\}
\end{array}$
\end{description}}
\vbox{\begin{description}
\item[\T{10}{64}:] [(1.42)]\\

$\begin{array}{@{}*{5}{l@{\hspace{4pt}}}l@{}}
 \{1,2,3\}& \{1,2,3\}& \{1,4,5\}& \{1,5,6\}& \{1,6,7\}& \{1,7,8\}\\
 \{1,7,8\}& \{1,8,9\}& \{0,1,9\}& \{0,1,4\}& \{2,4,8\}& \{2,4,5\}\\
 \{2,4,5\}& \{2,6,9\}& \{0,2,7\}& \{2,5,8\}& \{2,7,9\}& \{0,2,6\}\\
 \{0,2,6\}& \{3,4,7\}& \{3,5,9\}& \{0,3,6\}& \{3,5,7\}& \{3,6,8\}\\
 \{3,6,8\}& \{3,4,9\}& \{0,3,8\}& \{4,6,8\}& \{4,6,9\}& \{0,4,7\}
\end{array}$
\end{description}}
\vbox{\begin{description}
\item[\T{10}{65}:] [(1.12)]\\

$\begin{array}{@{}*{5}{l@{\hspace{4pt}}}l@{}}
 \{1,2,3\}& \{1,2,3\}& \{1,4,5\}& \{1,5,6\}& \{1,6,7\}& \{1,7,8\}\\
 \{1,7,8\}& \{1,8,9\}& \{0,1,9\}& \{0,1,4\}& \{2,4,6\}& \{2,5,7\}\\
 \{2,5,7\}& \{2,6,8\}& \{2,7,9\}& \{0,2,8\}& \{2,4,9\}& \{0,2,5\}\\
 \{0,2,5\}& \{3,4,8\}& \{3,5,9\}& \{0,3,6\}& \{3,4,7\}& \{3,5,8\}\\
 \{3,5,8\}& \{3,6,9\}& \{0,3,7\}& \{4,6,7\}& \{4,5,9\}& \{0,4,8\}
\end{array}$
\end{description}}
\vbox{\begin{description}
\item[\T{10}{66}:] [(1.44)]\\

$\begin{array}{@{}*{5}{l@{\hspace{4pt}}}l@{}}
 \{1,2,3\}& \{1,2,3\}& \{1,4,5\}& \{1,5,6\}& \{1,4,6\}& \{1,7,8\}\\
 \{1,7,8\}& \{1,8,9\}& \{0,1,9\}& \{0,1,7\}& \{2,4,7\}& \{2,5,8\}\\
 \{2,5,8\}& \{2,5,6\}& \{2,7,9\}& \{0,2,8\}& \{2,4,9\}& \{0,2,6\}\\
 \{0,2,6\}& \{0,3,4\}& \{3,5,9\}& \{3,6,7\}& \{3,5,7\}& \{3,4,8\}\\
 \{3,4,8\}& \{3,6,9\}& \{0,3,8\}& \{0,4,6\}& \{4,5,8\}& \{4,7,9\}
\end{array}$
\end{description}}
\vbox{\begin{description}
\item[\T{10}{67}:] [(1.25)]\\

$\begin{array}{@{}*{5}{l@{\hspace{4pt}}}l@{}}
 \{1,2,3\}& \{1,2,3\}& \{1,4,5\}& \{1,5,6\}& \{1,4,6\}& \{1,7,8\}\\
 \{1,7,8\}& \{1,8,9\}& \{0,1,9\}& \{0,1,7\}& \{2,4,6\}& \{2,5,8\}\\
 \{2,5,8\}& \{2,6,7\}& \{2,7,9\}& \{0,2,8\}& \{2,5,9\}& \{0,2,4\}\\
 \{0,2,4\}& \{3,4,7\}& \{0,3,5\}& \{3,6,9\}& \{3,5,7\}& \{3,4,8\}\\
 \{3,4,8\}& \{3,8,9\}& \{0,3,6\}& \{4,5,8\}& \{4,7,9\}& \{0,4,9\}
\end{array}$
\end{description}}
\vbox{\begin{description}
\item[\T{10}{68}:] [(1.3)]\\

$\begin{array}{@{}*{5}{l@{\hspace{4pt}}}l@{}}
 \{1,2,3\}& \{1,2,3\}& \{1,4,5\}& \{1,5,6\}& \{1,6,7\}& \{1,7,8\}\\
 \{1,7,8\}& \{1,8,9\}& \{0,1,9\}& \{0,1,4\}& \{2,4,6\}& \{2,5,8\}\\
 \{2,5,8\}& \{2,5,6\}& \{2,7,9\}& \{0,2,8\}& \{2,4,9\}& \{0,2,7\}\\
 \{0,2,7\}& \{3,4,9\}& \{3,5,7\}& \{3,6,8\}& \{0,3,7\}& \{3,4,8\}\\
 \{3,4,8\}& \{3,5,9\}& \{0,3,6\}& \{4,6,7\}& \{4,7,8\}& \{0,4,5\}
\end{array}$
\end{description}}
\vbox{\begin{description}
\item[\T{10}{69}:] [(1.20)]\\

$\begin{array}{@{}*{5}{l@{\hspace{4pt}}}l@{}}
 \{1,2,3\}& \{1,2,3\}& \{1,4,5\}& \{1,5,6\}& \{1,4,6\}& \{1,7,8\}\\
 \{1,7,8\}& \{1,8,9\}& \{0,1,9\}& \{0,1,7\}& \{2,4,6\}& \{2,5,9\}\\
 \{2,5,9\}& \{2,6,7\}& \{2,5,7\}& \{0,2,8\}& \{2,8,9\}& \{0,2,4\}\\
 \{0,2,4\}& \{3,4,8\}& \{0,3,5\}& \{3,6,9\}& \{3,4,7\}& \{3,5,8\}\\
 \{3,5,8\}& \{3,7,9\}& \{0,3,6\}& \{4,7,9\}& \{4,5,9\}& \{0,4,8\}
\end{array}$
\end{description}}
\vbox{\begin{description}
\item[\T{10}{70}:] [(1.46)]\\

$\begin{array}{@{}*{5}{l@{\hspace{4pt}}}l@{}}
 \{1,2,3\}& \{1,2,3\}& \{1,4,5\}& \{1,5,6\}& \{1,4,6\}& \{1,7,8\}\\
 \{1,7,8\}& \{1,8,9\}& \{0,1,9\}& \{0,1,7\}& \{2,4,7\}& \{0,2,5\}\\
 \{0,2,5\}& \{2,6,8\}& \{2,7,9\}& \{2,5,8\}& \{2,4,9\}& \{0,2,6\}\\
 \{0,2,6\}& \{3,4,8\}& \{3,5,7\}& \{3,6,9\}& \{3,6,7\}& \{0,3,8\}\\
 \{0,3,8\}& \{3,5,9\}& \{0,3,4\}& \{4,6,7\}& \{4,8,9\}& \{0,4,5\}
\end{array}$
\end{description}}
\vbox{\begin{description}
\item[\T{10}{71}:] [(1.35)]\\

$\begin{array}{@{}*{5}{l@{\hspace{4pt}}}l@{}}
 \{1,2,3\}& \{1,2,3\}& \{1,4,5\}& \{1,5,6\}& \{1,4,6\}& \{1,7,8\}\\
 \{1,7,8\}& \{1,8,9\}& \{0,1,9\}& \{0,1,7\}& \{2,4,7\}& \{2,5,8\}\\
 \{2,5,8\}& \{2,6,9\}& \{2,5,7\}& \{0,2,8\}& \{2,4,9\}& \{0,2,6\}\\
 \{0,2,6\}& \{3,4,8\}& \{0,3,5\}& \{3,6,7\}& \{3,7,9\}& \{3,6,8\}\\
 \{3,6,8\}& \{3,5,9\}& \{0,3,4\}& \{4,6,8\}& \{4,7,9\}& \{0,4,5\}
\end{array}$
\end{description}}
\vbox{\begin{description}
\item[\T{10}{72}:] [(1.14)]\\

$\begin{array}{@{}*{5}{l@{\hspace{4pt}}}l@{}}
 \{1,2,3\}& \{1,2,3\}& \{1,4,5\}& \{1,5,6\}& \{1,6,7\}& \{1,7,8\}\\
 \{1,7,8\}& \{1,8,9\}& \{0,1,9\}& \{0,1,4\}& \{2,4,6\}& \{2,5,8\}\\
 \{2,5,8\}& \{2,6,9\}& \{2,4,7\}& \{0,2,8\}& \{2,5,9\}& \{0,2,7\}\\
 \{0,2,7\}& \{3,4,9\}& \{0,3,5\}& \{3,6,8\}& \{3,5,7\}& \{3,4,8\}\\
 \{3,4,8\}& \{3,7,9\}& \{0,3,6\}& \{4,6,7\}& \{4,5,8\}& \{0,4,9\}
\end{array}$
\end{description}}
\vbox{\begin{description}
\item[\T{10}{73}:] [(1.43)]\\

$\begin{array}{@{}*{5}{l@{\hspace{4pt}}}l@{}}
 \{1,2,3\}& \{1,2,3\}& \{1,4,5\}& \{1,5,6\}& \{1,6,7\}& \{1,7,8\}\\
 \{1,7,8\}& \{1,8,9\}& \{0,1,9\}& \{0,1,4\}& \{2,4,8\}& \{2,5,7\}\\
 \{2,5,7\}& \{2,4,6\}& \{0,2,7\}& \{2,6,8\}& \{2,5,9\}& \{0,2,9\}\\
 \{0,2,9\}& \{3,4,9\}& \{0,3,5\}& \{3,6,8\}& \{3,4,7\}& \{3,5,8\}\\
 \{3,5,8\}& \{3,7,9\}& \{0,3,6\}& \{4,6,9\}& \{4,5,7\}& \{0,4,8\}
\end{array}$
\end{description}}
\vbox{\begin{description}
\item[\T{10}{74}:] [(1.26)]\\

$\begin{array}{@{}*{5}{l@{\hspace{4pt}}}l@{}}
 \{1,2,3\}& \{1,2,3\}& \{1,4,5\}& \{1,5,6\}& \{1,4,6\}& \{1,7,8\}\\
 \{1,7,8\}& \{1,8,9\}& \{0,1,9\}& \{0,1,7\}& \{2,4,6\}& \{0,2,5\}\\
 \{0,2,5\}& \{2,6,7\}& \{2,7,9\}& \{2,4,8\}& \{2,5,9\}& \{0,2,8\}\\
 \{0,2,8\}& \{3,4,7\}& \{3,5,9\}& \{3,6,8\}& \{0,3,7\}& \{3,5,8\}\\
 \{3,5,8\}& \{3,4,9\}& \{0,3,6\}& \{0,4,8\}& \{4,7,9\}& \{0,4,5\}
\end{array}$
\end{description}}
\vbox{\begin{description}
\item[\T{10}{75}:] [(1.28)]\\

$\begin{array}{@{}*{5}{l@{\hspace{4pt}}}l@{}}
 \{1,2,3\}& \{1,2,3\}& \{1,4,5\}& \{1,5,6\}& \{1,4,6\}& \{1,7,8\}\\
 \{1,7,8\}& \{1,8,9\}& \{0,1,9\}& \{0,1,7\}& \{2,4,6\}& \{0,2,5\}\\
 \{0,2,5\}& \{2,6,7\}& \{2,7,9\}& \{2,4,8\}& \{2,5,9\}& \{0,2,8\}\\
 \{0,2,8\}& \{3,4,9\}& \{3,5,8\}& \{0,3,6\}& \{3,4,7\}& \{3,7,8\}\\
 \{3,7,8\}& \{3,6,9\}& \{0,3,5\}& \{0,4,7\}& \{4,5,8\}& \{0,4,9\}
\end{array}$
\end{description}}
\vbox{\begin{description}
\item[\T{10}{76}:] [(1.37)]\\

$\begin{array}{@{}*{5}{l@{\hspace{4pt}}}l@{}}
 \{1,2,3\}& \{1,2,3\}& \{1,4,5\}& \{1,5,6\}& \{1,4,6\}& \{1,7,8\}\\
 \{1,7,8\}& \{1,8,9\}& \{0,1,9\}& \{0,1,7\}& \{2,4,7\}& \{2,5,9\}\\
 \{2,5,9\}& \{2,6,8\}& \{2,6,7\}& \{0,2,8\}& \{2,4,9\}& \{0,2,5\}\\
 \{0,2,5\}& \{0,3,4\}& \{3,5,7\}& \{3,6,9\}& \{3,4,7\}& \{3,6,8\}\\
 \{3,6,8\}& \{3,5,9\}& \{0,3,8\}& \{0,4,6\}& \{4,5,8\}& \{4,8,9\}
\end{array}$
\end{description}}
\vbox{\begin{description}
\item[\T{10}{77}:] [(1.8)]\\

$\begin{array}{@{}*{5}{l@{\hspace{4pt}}}l@{}}
 \{1,2,3\}& \{1,2,3\}& \{1,4,5\}& \{1,5,6\}& \{1,6,7\}& \{1,7,8\}\\
 \{1,7,8\}& \{1,8,9\}& \{0,1,9\}& \{0,1,4\}& \{2,4,6\}& \{0,2,5\}\\
 \{0,2,5\}& \{2,6,8\}& \{2,7,9\}& \{2,5,8\}& \{2,4,9\}& \{0,2,7\}\\
 \{0,2,7\}& \{3,4,7\}& \{3,5,9\}& \{0,3,6\}& \{3,5,7\}& \{3,4,8\}\\
 \{3,4,8\}& \{3,6,9\}& \{0,3,8\}& \{4,7,8\}& \{4,5,9\}& \{0,4,6\}
\end{array}$
\end{description}}
\vbox{\begin{description}
\item[\T{10}{78}:] [(1.11)]\\

$\begin{array}{@{}*{5}{l@{\hspace{4pt}}}l@{}}
 \{1,2,3\}& \{1,2,3\}& \{1,4,5\}& \{1,5,6\}& \{1,6,7\}& \{1,7,8\}\\
 \{1,7,8\}& \{1,8,9\}& \{0,1,9\}& \{0,1,4\}& \{2,4,6\}& \{2,5,7\}\\
 \{2,5,7\}& \{2,6,8\}& \{2,7,9\}& \{0,2,8\}& \{2,4,9\}& \{0,2,5\}\\
 \{0,2,5\}& \{3,4,7\}& \{3,5,8\}& \{0,3,6\}& \{3,6,7\}& \{3,4,8\}\\
 \{3,4,8\}& \{3,5,9\}& \{0,3,9\}& \{4,5,8\}& \{4,6,9\}& \{0,4,7\}
\end{array}$
\end{description}}
\vbox{\begin{description}
\item[\T{10}{79}:] [(2.21)]\\

$\begin{array}{@{}*{5}{l@{\hspace{4pt}}}l@{}}
 \{1,2,3\}& \{1,2,3\}& \{1,4,5\}& \{1,4,5\}& \{1,6,7\}& \{1,7,8\}\\
 \{1,7,8\}& \{1,8,9\}& \{0,1,9\}& \{0,1,6\}& \{2,4,6\}& \{2,5,8\}\\
 \{2,5,8\}& \{2,5,6\}& \{0,2,7\}& \{2,7,8\}& \{2,4,9\}& \{0,2,9\}\\
 \{0,2,9\}& \{3,4,7\}& \{3,5,9\}& \{3,6,8\}& \{3,5,7\}& \{0,3,8\}\\
 \{0,3,8\}& \{3,6,9\}& \{0,3,4\}& \{4,6,8\}& \{4,8,9\}& \{0,4,7\}
\end{array}$
\end{description}}
\vbox{\begin{description}
\item[\T{10}{80}:] [(2.23)]\\

$\begin{array}{@{}*{5}{l@{\hspace{4pt}}}l@{}}
 \{1,2,3\}& \{1,2,3\}& \{1,4,5\}& \{1,4,5\}& \{1,6,7\}& \{1,7,8\}\\
 \{1,7,8\}& \{1,8,9\}& \{0,1,9\}& \{0,1,6\}& \{2,4,6\}& \{2,5,8\}\\
 \{2,5,8\}& \{2,5,6\}& \{2,7,9\}& \{0,2,8\}& \{2,4,9\}& \{0,2,7\}\\
 \{0,2,7\}& \{3,4,8\}& \{3,5,9\}& \{0,3,6\}& \{3,5,7\}& \{3,7,8\}\\
 \{3,7,8\}& \{3,6,9\}& \{0,3,4\}& \{4,6,7\}& \{4,7,9\}& \{0,4,8\}
\end{array}$
\end{description}}
\vbox{\begin{description}
\item[\T{10}{81}:] [(2.27)]\\

$\begin{array}{@{}*{5}{l@{\hspace{4pt}}}l@{}}
 \{1,2,3\}& \{1,2,3\}& \{1,4,5\}& \{1,4,5\}& \{1,6,7\}& \{1,7,8\}\\
 \{1,7,8\}& \{1,8,9\}& \{0,1,9\}& \{0,1,6\}& \{2,4,6\}& \{2,5,9\}\\
 \{2,5,9\}& \{0,2,6\}& \{2,5,7\}& \{2,7,8\}& \{2,8,9\}& \{0,2,4\}\\
 \{0,2,4\}& \{3,4,9\}& \{3,5,8\}& \{3,5,6\}& \{3,4,7\}& \{0,3,8\}\\
 \{0,3,8\}& \{3,6,9\}& \{0,3,7\}& \{4,7,9\}& \{4,6,8\}& \{0,4,8\}
\end{array}$
\end{description}}
\vbox{\begin{description}
\item[\T{10}{82}:] [(2.24)]\\

$\begin{array}{@{}*{5}{l@{\hspace{4pt}}}l@{}}
 \{1,2,3\}& \{1,2,3\}& \{1,4,5\}& \{1,4,5\}& \{1,6,7\}& \{1,7,8\}\\
 \{1,7,8\}& \{1,8,9\}& \{0,1,9\}& \{0,1,6\}& \{2,4,6\}& \{2,5,7\}\\
 \{2,5,7\}& \{2,6,8\}& \{0,2,7\}& \{2,5,8\}& \{2,4,9\}& \{0,2,9\}\\
 \{0,2,9\}& \{3,4,8\}& \{3,5,6\}& \{0,3,6\}& \{3,7,9\}& \{3,7,8\}\\
 \{3,7,8\}& \{3,5,9\}& \{0,3,4\}& \{4,6,7\}& \{4,7,9\}& \{0,4,8\}
\end{array}$
\end{description}}
\vbox{\begin{description}
\item[\T{10}{83}:] [(2.10)]\\

$\begin{array}{@{}*{5}{l@{\hspace{4pt}}}l@{}}
 \{1,2,3\}& \{1,2,3\}& \{1,4,5\}& \{1,4,5\}& \{1,6,7\}& \{1,7,8\}\\
 \{1,7,8\}& \{1,8,9\}& \{0,1,9\}& \{0,1,6\}& \{2,4,6\}& \{2,5,8\}\\
 \{2,5,8\}& \{2,5,6\}& \{2,4,7\}& \{0,2,8\}& \{2,7,9\}& \{0,2,9\}\\
 \{0,2,9\}& \{3,4,8\}& \{3,5,9\}& \{0,3,6\}& \{3,6,7\}& \{3,7,8\}\\
 \{3,7,8\}& \{3,4,9\}& \{0,3,5\}& \{0,4,7\}& \{4,6,9\}& \{0,4,8\}
\end{array}$
\end{description}}
\vbox{\begin{description}
\item[\T{10}{84}:] [(2.14)]\\

$\begin{array}{@{}*{5}{l@{\hspace{4pt}}}l@{}}
 \{1,2,3\}& \{1,2,3\}& \{1,4,5\}& \{1,4,5\}& \{1,6,7\}& \{1,7,8\}\\
 \{1,7,8\}& \{1,8,9\}& \{0,1,9\}& \{0,1,6\}& \{2,4,6\}& \{2,5,7\}\\
 \{2,5,7\}& \{2,6,8\}& \{0,2,7\}& \{2,5,8\}& \{2,4,9\}& \{0,2,9\}\\
 \{0,2,9\}& \{3,4,7\}& \{0,3,5\}& \{3,4,6\}& \{3,7,9\}& \{3,6,8\}\\
 \{3,6,8\}& \{3,5,9\}& \{0,3,8\}& \{0,4,8\}& \{4,8,9\}& \{0,4,7\}
\end{array}$
\end{description}}
\vbox{\begin{description}
\item[\T{10}{85}:] [(2.12)]\\

$\begin{array}{@{}*{5}{l@{\hspace{4pt}}}l@{}}
 \{1,2,3\}& \{1,2,3\}& \{1,4,5\}& \{1,4,5\}& \{1,6,7\}& \{1,7,8\}\\
 \{1,7,8\}& \{1,8,9\}& \{0,1,9\}& \{0,1,6\}& \{2,4,6\}& \{2,5,9\}\\
 \{2,5,9\}& \{2,5,6\}& \{2,4,7\}& \{0,2,8\}& \{2,8,9\}& \{0,2,7\}\\
 \{0,2,7\}& \{3,4,9\}& \{3,5,7\}& \{3,6,8\}& \{0,3,7\}& \{3,5,8\}\\
 \{3,5,8\}& \{3,6,9\}& \{0,3,4\}& \{4,7,9\}& \{4,6,8\}& \{0,4,8\}
\end{array}$
\end{description}}
\vbox{\begin{description}
\item[\T{10}{86}:] [(2.4)]\\

$\begin{array}{@{}*{5}{l@{\hspace{4pt}}}l@{}}
 \{1,2,3\}& \{1,2,3\}& \{1,4,5\}& \{1,4,5\}& \{1,6,7\}& \{1,7,8\}\\
 \{1,7,8\}& \{1,8,9\}& \{0,1,9\}& \{0,1,6\}& \{2,4,6\}& \{2,5,9\}\\
 \{2,5,9\}& \{2,5,6\}& \{0,2,7\}& \{2,4,8\}& \{2,7,9\}& \{0,2,8\}\\
 \{0,2,8\}& \{3,4,7\}& \{0,3,5\}& \{3,6,8\}& \{3,6,7\}& \{3,5,8\}\\
 \{3,5,8\}& \{3,4,9\}& \{0,3,9\}& \{0,4,8\}& \{4,6,9\}& \{0,4,7\}
\end{array}$
\end{description}}
\vbox{\begin{description}
\item[\T{10}{87}:] [(2.6)]\\

$\begin{array}{@{}*{5}{l@{\hspace{4pt}}}l@{}}
 \{1,2,3\}& \{1,2,3\}& \{1,4,5\}& \{1,4,5\}& \{1,6,7\}& \{1,7,8\}\\
 \{1,7,8\}& \{1,8,9\}& \{0,1,9\}& \{0,1,6\}& \{2,4,6\}& \{2,5,9\}\\
 \{2,5,9\}& \{2,5,6\}& \{2,4,7\}& \{0,2,8\}& \{2,8,9\}& \{0,2,7\}\\
 \{0,2,7\}& \{3,4,8\}& \{3,5,7\}& \{0,3,6\}& \{3,7,9\}& \{3,6,8\}\\
 \{3,6,8\}& \{3,4,9\}& \{0,3,5\}& \{0,4,7\}& \{4,6,9\}& \{0,4,8\}
\end{array}$
\end{description}}
\vbox{\begin{description}
\item[\T{10}{88}:] [(2.22)]\\

$\begin{array}{@{}*{5}{l@{\hspace{4pt}}}l@{}}
 \{1,2,3\}& \{1,2,3\}& \{1,4,5\}& \{1,4,5\}& \{1,6,7\}& \{1,7,8\}\\
 \{1,7,8\}& \{1,8,9\}& \{0,1,9\}& \{0,1,6\}& \{2,4,6\}& \{2,5,8\}\\
 \{2,5,8\}& \{2,5,6\}& \{0,2,7\}& \{2,7,8\}& \{2,4,9\}& \{0,2,9\}\\
 \{0,2,9\}& \{0,3,4\}& \{3,5,9\}& \{3,6,8\}& \{3,5,7\}& \{3,4,8\}\\
 \{3,4,8\}& \{3,6,9\}& \{0,3,7\}& \{4,6,7\}& \{0,4,8\}& \{4,7,9\}
\end{array}$
\end{description}}
\vbox{\begin{description}
\item[\T{10}{89}:] [(2.5)]\\

$\begin{array}{@{}*{5}{l@{\hspace{4pt}}}l@{}}
 \{1,2,3\}& \{1,2,3\}& \{1,4,5\}& \{1,4,5\}& \{1,6,7\}& \{1,7,8\}\\
 \{1,7,8\}& \{1,8,9\}& \{0,1,9\}& \{0,1,6\}& \{2,4,6\}& \{2,5,9\}\\
 \{2,5,9\}& \{2,5,6\}& \{0,2,7\}& \{2,4,8\}& \{2,7,9\}& \{0,2,8\}\\
 \{0,2,8\}& \{3,4,9\}& \{0,3,5\}& \{3,6,8\}& \{3,4,7\}& \{3,5,8\}\\
 \{3,5,8\}& \{3,6,9\}& \{0,3,7\}& \{4,6,7\}& \{0,4,8\}& \{0,4,9\}
\end{array}$
\end{description}}
\vbox{\begin{description}
\item[\T{10}{90}:] [(2.3)]\\

$\begin{array}{@{}*{5}{l@{\hspace{4pt}}}l@{}}
 \{1,2,3\}& \{1,2,3\}& \{1,4,5\}& \{1,4,5\}& \{1,6,7\}& \{1,7,8\}\\
 \{1,7,8\}& \{1,8,9\}& \{0,1,9\}& \{0,1,6\}& \{2,4,6\}& \{2,5,7\}\\
 \{2,5,7\}& \{2,5,6\}& \{0,2,7\}& \{2,4,8\}& \{2,8,9\}& \{0,2,9\}\\
 \{0,2,9\}& \{3,4,8\}& \{0,3,5\}& \{3,6,9\}& \{3,6,7\}& \{3,7,8\}\\
 \{3,7,8\}& \{3,5,9\}& \{0,3,4\}& \{4,7,9\}& \{4,6,9\}& \{0,4,7\}
\end{array}$
\end{description}}
\vbox{\begin{description}
\item[\T{10}{91}:] [(2.30)]\\

$\begin{array}{@{}*{5}{l@{\hspace{4pt}}}l@{}}
 \{1,2,3\}& \{1,2,3\}& \{1,4,5\}& \{1,5,6\}& \{1,6,7\}& \{1,7,8\}\\
 \{1,7,8\}& \{1,8,9\}& \{0,1,9\}& \{0,1,4\}& \{2,4,6\}& \{0,2,5\}\\
 \{0,2,5\}& \{2,5,6\}& \{2,4,7\}& \{2,7,8\}& \{2,8,9\}& \{0,2,9\}\\
 \{0,2,9\}& \{3,4,7\}& \{3,5,8\}& \{3,6,9\}& \{0,3,7\}& \{3,6,8\}\\
 \{3,6,8\}& \{3,4,9\}& \{0,3,5\}& \{4,5,8\}& \{4,6,9\}& \{0,4,8\}
\end{array}$
\end{description}}
\vbox{\begin{description}
\item[\T{10}{92}:] [(2.19)]\\

$\begin{array}{@{}*{5}{l@{\hspace{4pt}}}l@{}}
 \{1,2,3\}& \{1,2,3\}& \{1,4,5\}& \{1,4,5\}& \{1,6,7\}& \{1,7,8\}\\
 \{1,7,8\}& \{1,8,9\}& \{0,1,9\}& \{0,1,6\}& \{2,4,6\}& \{0,2,5\}\\
 \{0,2,5\}& \{2,6,8\}& \{2,7,9\}& \{2,5,8\}& \{2,4,9\}& \{0,2,7\}\\
 \{0,2,7\}& \{3,4,9\}& \{3,5,8\}& \{0,3,6\}& \{3,4,7\}& \{3,7,8\}\\
 \{3,7,8\}& \{3,6,9\}& \{0,3,5\}& \{0,4,7\}& \{4,6,8\}& \{0,4,8\}
\end{array}$
\end{description}}
\vbox{\begin{description}
\item[\T{10}{93}:] [(2.18)]\\

$\begin{array}{@{}*{5}{l@{\hspace{4pt}}}l@{}}
 \{1,2,3\}& \{1,2,3\}& \{1,4,5\}& \{1,4,5\}& \{1,6,7\}& \{1,7,8\}\\
 \{1,7,8\}& \{1,8,9\}& \{0,1,9\}& \{0,1,6\}& \{2,4,6\}& \{0,2,5\}\\
 \{0,2,5\}& \{2,6,8\}& \{2,7,9\}& \{2,5,8\}& \{2,4,9\}& \{0,2,7\}\\
 \{0,2,7\}& \{3,4,8\}& \{3,5,9\}& \{0,3,6\}& \{3,4,7\}& \{3,6,8\}\\
 \{3,6,8\}& \{3,7,9\}& \{0,3,5\}& \{0,4,7\}& \{4,6,9\}& \{0,4,8\}
\end{array}$
\end{description}}
\vbox{\begin{description}
\item[\T{10}{94}:] [(2.26)]\\

$\begin{array}{@{}*{5}{l@{\hspace{4pt}}}l@{}}
 \{1,2,3\}& \{1,2,3\}& \{1,4,5\}& \{1,4,5\}& \{1,6,7\}& \{1,7,8\}\\
 \{1,7,8\}& \{1,8,9\}& \{0,1,9\}& \{0,1,6\}& \{2,4,6\}& \{2,5,7\}\\
 \{2,5,7\}& \{2,6,8\}& \{0,2,7\}& \{2,4,8\}& \{2,5,9\}& \{0,2,9\}\\
 \{0,2,9\}& \{0,3,4\}& \{3,5,6\}& \{3,6,9\}& \{3,4,7\}& \{3,7,8\}\\
 \{3,7,8\}& \{3,8,9\}& \{0,3,5\}& \{4,7,9\}& \{0,4,8\}& \{4,6,9\}
\end{array}$
\end{description}}
\vbox{\begin{description}
\item[\T{10}{95}:] [(2.9)]\\

$\begin{array}{@{}*{5}{l@{\hspace{4pt}}}l@{}}
 \{1,2,3\}& \{1,2,3\}& \{1,4,5\}& \{1,4,5\}& \{1,6,7\}& \{1,7,8\}\\
 \{1,7,8\}& \{1,8,9\}& \{0,1,9\}& \{0,1,6\}& \{2,4,6\}& \{2,5,8\}\\
 \{2,5,8\}& \{2,5,6\}& \{2,4,7\}& \{0,2,8\}& \{2,7,9\}& \{0,2,9\}\\
 \{0,2,9\}& \{3,4,8\}& \{3,5,7\}& \{0,3,6\}& \{3,7,9\}& \{3,6,8\}\\
 \{3,6,8\}& \{3,4,9\}& \{0,3,5\}& \{0,4,7\}& \{4,6,9\}& \{0,4,8\}
\end{array}$
\end{description}}
\vbox{\begin{description}
\item[\T{10}{96}:] [(2.8)]\\

$\begin{array}{@{}*{5}{l@{\hspace{4pt}}}l@{}}
 \{1,2,3\}& \{1,2,3\}& \{1,4,5\}& \{1,4,5\}& \{1,6,7\}& \{1,7,8\}\\
 \{1,7,8\}& \{1,8,9\}& \{0,1,9\}& \{0,1,6\}& \{2,4,6\}& \{2,5,8\}\\
 \{2,5,8\}& \{2,5,6\}& \{2,4,7\}& \{0,2,8\}& \{2,7,9\}& \{0,2,9\}\\
 \{0,2,9\}& \{3,4,8\}& \{3,5,7\}& \{3,6,9\}& \{0,3,7\}& \{3,6,8\}\\
 \{3,6,8\}& \{3,5,9\}& \{0,3,4\}& \{4,7,9\}& \{4,6,9\}& \{0,4,8\}
\end{array}$
\end{description}}
\vbox{\begin{description}
\item[\T{10}{97}:] [(2.15)]\\

$\begin{array}{@{}*{5}{l@{\hspace{4pt}}}l@{}}
 \{1,2,3\}& \{1,2,3\}& \{1,4,5\}& \{1,4,5\}& \{1,6,7\}& \{1,7,8\}\\
 \{1,7,8\}& \{1,8,9\}& \{0,1,9\}& \{0,1,6\}& \{2,4,6\}& \{2,5,9\}\\
 \{2,5,9\}& \{2,6,8\}& \{0,2,7\}& \{2,5,8\}& \{2,7,9\}& \{0,2,4\}\\
 \{0,2,4\}& \{3,4,9\}& \{3,5,8\}& \{0,3,6\}& \{3,4,7\}& \{3,7,8\}\\
 \{3,7,8\}& \{3,6,9\}& \{0,3,5\}& \{4,7,9\}& \{4,6,8\}& \{0,4,8\}
\end{array}$
\end{description}}
\vbox{\begin{description}
\item[\T{10}{98}:] [(2.11)]\\

$\begin{array}{@{}*{5}{l@{\hspace{4pt}}}l@{}}
 \{1,2,3\}& \{1,2,3\}& \{1,4,5\}& \{1,4,5\}& \{1,6,7\}& \{1,7,8\}\\
 \{1,7,8\}& \{1,8,9\}& \{0,1,9\}& \{0,1,6\}& \{2,4,6\}& \{2,5,8\}\\
 \{2,5,8\}& \{2,5,6\}& \{2,4,7\}& \{0,2,8\}& \{2,7,9\}& \{0,2,9\}\\
 \{0,2,9\}& \{3,4,9\}& \{0,3,5\}& \{3,6,8\}& \{3,5,7\}& \{3,7,8\}\\
 \{3,7,8\}& \{3,6,9\}& \{0,3,4\}& \{4,7,9\}& \{4,6,8\}& \{0,4,8\}
\end{array}$
\end{description}}
\vbox{\begin{description}
\item[\T{10}{99}:] [(2.20)]\\

$\begin{array}{@{}*{5}{l@{\hspace{4pt}}}l@{}}
 \{1,2,3\}& \{1,2,3\}& \{1,4,5\}& \{1,4,5\}& \{1,6,7\}& \{1,7,8\}\\
 \{1,7,8\}& \{1,8,9\}& \{0,1,9\}& \{0,1,6\}& \{2,4,6\}& \{0,2,5\}\\
 \{0,2,5\}& \{2,6,8\}& \{2,5,7\}& \{2,7,8\}& \{2,4,9\}& \{0,2,9\}\\
 \{0,2,9\}& \{0,3,4\}& \{3,5,6\}& \{3,6,9\}& \{3,4,7\}& \{3,5,8\}\\
 \{3,5,8\}& \{3,7,9\}& \{0,3,8\}& \{0,4,7\}& \{4,6,8\}& \{4,8,9\}
\end{array}$
\end{description}}
\vbox{\begin{description}
\item[\T{10}{100}:] [(2.17)]\\

$\begin{array}{@{}*{5}{l@{\hspace{4pt}}}l@{}}
 \{1,2,3\}& \{1,2,3\}& \{1,4,5\}& \{1,4,5\}& \{1,6,7\}& \{1,7,8\}\\
 \{1,7,8\}& \{1,8,9\}& \{0,1,9\}& \{0,1,6\}& \{2,4,6\}& \{0,2,5\}\\
 \{0,2,5\}& \{2,6,8\}& \{2,7,9\}& \{2,5,8\}& \{2,4,9\}& \{0,2,7\}\\
 \{0,2,7\}& \{3,4,8\}& \{3,5,7\}& \{3,6,9\}& \{0,3,7\}& \{3,6,8\}\\
 \{3,6,8\}& \{3,5,9\}& \{0,3,4\}& \{4,6,7\}& \{4,7,9\}& \{0,4,8\}
\end{array}$
\end{description}}
\vbox{\begin{description}
\item[\T{10}{101}:] [(2.16)]\\

$\begin{array}{@{}*{5}{l@{\hspace{4pt}}}l@{}}
 \{1,2,3\}& \{1,2,3\}& \{1,4,5\}& \{1,4,5\}& \{1,6,7\}& \{1,7,8\}\\
 \{1,7,8\}& \{1,8,9\}& \{0,1,9\}& \{0,1,6\}& \{2,4,6\}& \{0,2,5\}\\
 \{0,2,5\}& \{2,6,8\}& \{2,7,9\}& \{2,5,8\}& \{2,4,9\}& \{0,2,7\}\\
 \{0,2,7\}& \{3,4,7\}& \{3,5,9\}& \{3,5,6\}& \{3,6,7\}& \{0,3,8\}\\
 \{0,3,8\}& \{3,8,9\}& \{0,3,4\}& \{4,6,8\}& \{4,7,9\}& \{0,4,8\}
\end{array}$
\end{description}}
\vbox{\begin{description}
\item[\T{10}{102}:] [(2.28)]\\

$\begin{array}{@{}*{5}{l@{\hspace{4pt}}}l@{}}
 \{1,2,3\}& \{1,2,3\}& \{1,4,5\}& \{1,4,5\}& \{1,6,7\}& \{1,7,8\}\\
 \{1,7,8\}& \{1,8,9\}& \{0,1,9\}& \{0,1,6\}& \{2,4,6\}& \{2,5,8\}\\
 \{2,5,8\}& \{2,6,9\}& \{0,2,7\}& \{2,7,8\}& \{2,4,9\}& \{0,2,5\}\\
 \{0,2,5\}& \{3,4,7\}& \{3,5,6\}& \{3,6,8\}& \{3,7,9\}& \{0,3,8\}\\
 \{0,3,8\}& \{3,5,9\}& \{0,3,4\}& \{4,6,8\}& \{4,8,9\}& \{0,4,7\}
\end{array}$
\end{description}}
\vbox{\begin{description}
\item[\T{10}{103}:] [(2.29)]\\

$\begin{array}{@{}*{5}{l@{\hspace{4pt}}}l@{}}
 \{1,2,3\}& \{1,2,3\}& \{1,4,5\}& \{1,4,5\}& \{1,6,7\}& \{1,7,8\}\\
 \{1,7,8\}& \{1,8,9\}& \{0,1,9\}& \{0,1,6\}& \{2,4,6\}& \{0,2,5\}\\
 \{0,2,5\}& \{2,6,8\}& \{2,7,9\}& \{2,4,8\}& \{2,5,9\}& \{0,2,7\}\\
 \{0,2,7\}& \{0,3,4\}& \{3,5,7\}& \{3,6,9\}& \{3,6,7\}& \{3,5,8\}\\
 \{3,5,8\}& \{3,4,9\}& \{0,3,8\}& \{0,4,7\}& \{4,6,8\}& \{4,7,9\}
\end{array}$
\end{description}}
\vbox{\begin{description}
\item[\T{10}{104}:] [(2.13)]\\

$\begin{array}{@{}*{5}{l@{\hspace{4pt}}}l@{}}
 \{1,2,3\}& \{1,2,3\}& \{1,4,5\}& \{1,4,5\}& \{1,6,7\}& \{1,7,8\}\\
 \{1,7,8\}& \{1,8,9\}& \{0,1,9\}& \{0,1,6\}& \{2,4,6\}& \{2,5,7\}\\
 \{2,5,7\}& \{2,6,8\}& \{0,2,7\}& \{2,5,8\}& \{2,4,9\}& \{0,2,9\}\\
 \{0,2,9\}& \{3,4,7\}& \{0,3,5\}& \{3,6,9\}& \{3,6,7\}& \{3,4,8\}\\
 \{3,4,8\}& \{3,5,9\}& \{0,3,8\}& \{0,4,8\}& \{4,6,9\}& \{0,4,7\}
\end{array}$
\end{description}}
\vbox{\begin{description}
\item[\T{10}{105}:] [(2.7)]\\

$\begin{array}{@{}*{5}{l@{\hspace{4pt}}}l@{}}
 \{1,2,3\}& \{1,2,3\}& \{1,4,5\}& \{1,4,5\}& \{1,6,7\}& \{1,7,8\}\\
 \{1,7,8\}& \{1,8,9\}& \{0,1,9\}& \{0,1,6\}& \{2,4,6\}& \{0,2,5\}\\
 \{0,2,5\}& \{2,5,6\}& \{2,4,7\}& \{2,7,8\}& \{2,8,9\}& \{0,2,9\}\\
 \{0,2,9\}& \{3,4,9\}& \{3,5,7\}& \{3,6,8\}& \{0,3,7\}& \{3,5,8\}\\
 \{3,5,8\}& \{3,6,9\}& \{0,3,4\}& \{4,7,9\}& \{4,6,8\}& \{0,4,8\}
\end{array}$
\end{description}}
\vbox{\begin{description}
\item[\T{10}{106}:] [(2.25)]\\

$\begin{array}{@{}*{5}{l@{\hspace{4pt}}}l@{}}
 \{1,2,3\}& \{1,2,3\}& \{1,4,5\}& \{1,4,5\}& \{1,6,7\}& \{1,7,8\}\\
 \{1,7,8\}& \{1,8,9\}& \{0,1,9\}& \{0,1,6\}& \{2,4,6\}& \{0,2,5\}\\
 \{0,2,5\}& \{2,5,6\}& \{2,4,7\}& \{2,7,8\}& \{2,8,9\}& \{0,2,9\}\\
 \{0,2,9\}& \{3,4,8\}& \{3,5,7\}& \{3,6,9\}& \{0,3,7\}& \{3,6,8\}\\
 \{3,6,8\}& \{3,5,9\}& \{0,3,4\}& \{4,7,9\}& \{4,6,9\}& \{0,4,8\}
\end{array}%
$
\end{description}}
\subsection{The \texorpdfstring{$\TTS(10)$}{TTS(10)} and the good sequencings found.}\label{TabII10}

\vbox{\begin{description}
\item[\M{10}{1}{.1}:]\quad\\
$%
$

	Lexicographic least 3-good sequencing : 0123546789

	Number of 3-good sequencing found: 717200

	Lexicographic least 4-good sequencing : 0142756839

	Number of 4-good sequencing found: 100

	Number of 5-good sequencing found: 0
\end{description}}
\vbox{\begin{description}
\item[\M{10}{2}{.1}:]\quad\\
$%
$

	Lexicographic least 3-good sequencing : 0124356789

	Number of 3-good sequencing found: 715140

	Lexicographic least 4-good sequencing : 0152483967

	Number of 4-good sequencing found: 510

	Number of 5-good sequencing found: 0
\end{description}}
\vbox{\begin{description}
\item[\M{10}{3}{.1}:]\quad\\
$%
$

	Lexicographic least 3-good sequencing : 0124356789

	Number of 3-good sequencing found: 714200

	Lexicographic least 4-good sequencing : 0163248759

	Number of 4-good sequencing found: 620

	Number of 5-good sequencing found: 0
\end{description}}
\vbox{\begin{description}
\item[\M{10}{4}{.1}:]\quad\\
$%
$

	Lexicographic least 3-good sequencing : 0124356879

	Number of 3-good sequencing found: 710880

	Number of 4-good sequencing found: 0

	Number of 5-good sequencing found: 0
\end{description}}
\vbox{\begin{description}
\item[\M{10}{116}{.2}:]\quad\\
$%
$

	Lexicographic least 3-good sequencing : 0124356879

	Number of 3-good sequencing found: 710910

	Number of 4-good sequencing found: 0

	Number of 5-good sequencing found: 0
\end{description}}
\vbox{\begin{description}
\item[\M{10}{117}{.1}:]\quad\\
$%
$

	Lexicographic least 3-good sequencing : 0124356879

	Number of 3-good sequencing found: 708780

	Number of 4-good sequencing found: 0

	Number of 5-good sequencing found: 0
\end{description}}
\vbox{\begin{description}
\item[\M{10}{119}{.1}:]\quad\\
$%
$

	Lexicographic least 3-good sequencing : 0124356987

	Number of 3-good sequencing found: 691740

	Number of 4-good sequencing found: 0

	Number of 5-good sequencing found: 0
\end{description}}
\newpage
\subsection{The five special \texorpdfstring{$\MTS(10)$}{MTS(10)}.}\label{Tab29}
In this section we provide a 4-good sequencing for the  
``partial'' $\MTS(10)$ obtained by deleting a triple from
one of the five $\MTS(10)$, that do not have a 4-good sequencings.
\subsubsection{\texorpdfstring{$\M{10}{116}{.1}$}{M(10)116.1}}
%
\subsubsection{\texorpdfstring{$\M{10}{116}{.2}$}{M(10)116.2}}
%
\subsubsection{\texorpdfstring{$\M{10}{118}{.1}$}{M(10)118.1}}
%
\subsubsection{\texorpdfstring{$\M{10}{134}{.1}$}{M(10)134.1}}
%
\subsubsection{\texorpdfstring{$\M{10}{134}{.2}$}{M(10)134.2}}
%

\end{document}